\documentclass[10pt]{amsart}
\usepackage[T1]{fontenc}
\usepackage{geometry}
\usepackage[latin1] {inputenc}
\usepackage{amsmath}
\usepackage{amsfonts, amssymb, textcomp}
\usepackage[colorlinks=flase, linkcolor=red,urlcolor=green, citecolor=blue]{hyperref}
\usepackage{subeqnarray}

\usepackage{latexsym}
\usepackage{fancyhdr}
\usepackage{longtable}
\usepackage{amsmath, amssymb}
\usepackage{graphicx}
\usepackage{enumerate}

\numberwithin{equation}{section}

\theoremstyle{plain}
\newtheorem{proposition}{Proposition}[section]
\newtheorem{theorem}[proposition]{Theorem}
\newtheorem{lemma}[proposition]{Lemma}

\newtheorem{remark}[proposition]{Remark}

\newcommand{\RR}{\mathbb{R}}
\newcommand{\CC}{\mathbb{C}}
\newcommand{\NN}{\mathbb{N}}

\newcommand{\supp}{\operatorname{supp}}

\let\on=\operatorname

\newenvironment{demo}[1]{\par\smallskip\noindent{\bf #1.}}{\par\smallskip}

\makeatletter
\@namedef{subjclassname@2020}{%
  \textup{2020} Mathematics Subject Classification}
\makeatother

\newsavebox{\fmbox}
\newenvironment{fmpage}[1]
 {\begin{lrbox}{\fmbox}\begin{minipage}{#1}}
 {\end{minipage}\end{lrbox}\fbox{\usebox{\fmbox}}}

\title[On the maximal extension]
{On the maximal extension in the mixed ultradifferentiable weight sequence setting}

\author[G.~Schindl]{Gerhard Schindl}

\address{G.~Schindl: Fakult\"at f\"ur Mathematik, Universit\"at Wien, Oskar-Morgenstern-Platz~1, A-1090 Wien, Austria.}
\email{gerhard.schindl@univie.ac.at}

\begin{document}

\begin{abstract}
For the ultradifferentiable weight sequence setting it is known that the Borel map which assigns to each function the infinite jet of derivatives (at $0$) is surjective onto the corresponding weighted sequence class if and only if the sequence is strongly nonquasianalytic for both the Roumieu- and Beurling-type classes. Sequences which are nonquasianalytic but not strongly nonquasianalytic admit a controlled loss of regularity and we determine the maximal sequence for which such a mixed setting is possible for both types, hence get information on the controlled loss of surjectivity in this situation. Moreover, we compare the optimal sequences for both mixed strong nonquasianalyticity conditions arising in the literature.
\end{abstract}

\thanks{G. Schindl is supported by FWF-Projects P32905-N and P33417-N}
\date{\today}

\maketitle

\section{Introduction}\label{Introduction}
The study of the injectivity and surjectivity of the Borel map $\mathcal{B}: f\mapsto(f^{(j)}(0))_{j\in\NN}$ in the ultradifferentiable setting has a long tradition and these properties have been fully characterized in terms of the defining weight sequence $M=(M_j)_j\in\RR_{>0}^{\NN}$, assuming some standard growth and regularity properties (and similarly for weight functions $\omega$ as well). We recall some facts, for more detailed definitions and conditions see Section \ref{notation}. $\mathcal{B}$ is defined on classes $\mathcal{E}_{\{M\}}$ resp. $\mathcal{E}_{(M)}$ of {\itshape Roumieu-type} resp. {\itshape Beurling-type} and the target spaces of weighted sequences of (complex) numbers are denoted by $\Lambda_{\{M\}}$ resp. $\Lambda_{(M)}$. We denote by $[\cdot]$ ultradifferentiable classes either of Roumieu-type or of {\itshape Beurling-type}, but not mixing the cases, and similarly for the weighted sequence classes.

For our results we will assume the following basic properties for $M$: $1=M_0\le M_1$, log-convexity for $M$ and $\lim_{j\rightarrow+\infty}(M_j)^{1/j}=+\infty$ (for short we write $M\in\hyperlink{LCset}{\mathcal{LC}}$ for this class).\vspace{6pt}

For the injectivity of $\mathcal{B}$, characterized in the {\itshape Denjoy-Carleman Theorem}, we refer to \cite[Theorem 1.3.8]{hoermander} and to \cite[Theorem 4.2]{Komatsu73}. The ultradifferentiable class is quasianalytic, i.e. $\mathcal{B}$ defined on this class is injective, if and only if the defining sequence $M$ is quasianalytic, i.e. condition \hyperlink{mnq}{$(\text{nq})$} in Section \ref{mixedsnqs} is violated. Characterizing surjectivity has been treated in \cite{petzsche}, here the crucial condition is \hyperlink{gamma1}{$(\gamma_1)$} also known as ''strong nonquasianalyticity condition''. Both characterizations are valid for both the Roumieu- and Beurling-type classes and nonquasianalyticity characterizes the nontriviality for the test function space $\mathcal{D}_{[M]}([-1,1])$ (see again Section \ref{notation} for definition).

Moreover, it has turned out that in the quasianalytic setting where the function class strictly contains the real analytic functions, for both the Beurling- and Roumieu-type classes the Borel map can never be surjective onto the corresponding sequence classes (e.g. see \cite[Theorem 3]{thilliez} and \cite{borelmappingquasianalytic}). In \cite{Borelmapgenericity} and \cite{Borelmapalgebraity} it has been shown that the image of $\mathcal{B}$ is small and the ''size of the failure'' has been measured by using different concepts; for more details we refer to the introductions and citations in those papers. However, it is still an open question to give a full characterization of the image of $\mathcal{B}$ in this setting, i.e. deciding whether a given sequence of complex numbers belongs to the image or not.\vspace{6pt}

When the sequence $N$ is {\itshape nonquasianalytic but not strongly nonquasianalytic,} (in which case the real-analytic functions are always strictly contained in the function class), from the known results it follows that $\mathcal{B}$ is neither injective nor surjective (for both types). However, it is also known that for this case a controlled loss of regularity is possible (''mixed settings'') meaning that $\mathcal{B}(\mathcal{D}_{[N]}([-1,1]))\supseteq\Lambda_{[M]}$ for a different weight $M\le N$. In this context, when $N$ is fixed, it is natural to ask for the {\itshape maximal sequence $M$ admitting this inclusion} and this is the main problem studied in this article.

In \cite{surjectivity} a precise characterization of the inclusion $\mathcal{B}(\mathcal{D}_{[N]}([-1,1]))\supseteq\Lambda_{[M]}$, given in terms of two sequences $M$ and $N$, has been obtained. The crucial condition reads
$$\exists\;s\in\NN_{>0}:\;\;\sup_{p\in\NN_{>0}}\frac{\lambda_{p,s}^{M,N}}{p}\sum_{k\ge p}\frac{1}{\nu_k}<+\infty,$$
with $\lambda^{M,N}_{p,s}:=\sup_{0\le j<p}\left(\frac{M_p}{s^p N_j}\right)^{1/(p-j)}$ and $\nu_k:=\frac{N_k}{N_{k-1}}$. In the following this requirement is denoted by $(M,N)_{SV}$. In \cite{mixedramisurj} the results from \cite{surjectivity} have been generalized by involving a ramification parameter $r\in\NN_{>0}$ and dealing with special classes of ultradifferentiable ramification spaces introduced in \cite{Schmetsvaldivia00} and needed for the study of the surjectivity of the (asymptotic) Borel map in the ultraholomorphic setting, see Section \ref{ramisection}.

But in the literature there exists a second relevant and in general stronger mixed condition given in terms of $M$ and $N$, namely
$$\sup_{p\in\NN_{>0}}\frac{\mu_p}{p}\sum_{k\ge p}\frac{1}{\nu_k}<+\infty,$$
denoted by $(M,N)_{\gamma_1}$ in this work (with $\mu_j:=\frac{M_j}{M_{j-1}}$). $(M,N)_{\gamma_1}$ is easier to handle and a ''more natural'' generalization of \hyperlink{gamma1}{$(\gamma_1)$}. It has been introduced in \cite{ChaumatChollet94} and also used in \cite{whitneyextensionweightmatrix}. Finally, in \cite{mixedramisurj} and \cite{mixedsectorialextensions} a ramified generalization of this condition appeared, see again Section \ref{ramisection}.

Note that in $(M,N)_{\gamma_1}$ the sequence of quotients of $M$, denoted by $(\mu_j)_j$, appears whereas $(M,N)_{SV}$ is connected with the sequence of roots $(M_j)^{1/j}$ and these sequences are comparable (up to a constant) if and only if $M$ satisfies the technical assumption of {\itshape moderate growth} (by \cite[Lemma 2.2]{whitneyextensionweightmatrix}).\vspace{6pt}

By the construction from \cite[Sect. 4.1]{whitneyextensionweightmatrix} the optimal (i.e. largest) sequence $M\le N$ expressed in terms of $(M,N)_{\gamma_1}$, when $N$ is fixed, is already known; it is called the {\itshape descendant}, see Section \ref{descendant}. However, the main results of \cite{surjectivity} (and \cite{mixedramisurj}) show that $(M,N)_{SV}$ is the correct (but more technical) characterizing condition in the mixed setting, and the optimal sequence expressed in terms of this condition has not been computed so far; this is the main goal of this article. The main result in this context reads as follows, see Theorem \ref{maximalelement}:

\begin{theorem}\label{maximalelementintro}
Let $N\in\hyperlink{LCset}{\mathcal{LC}}$. Then the set of sequences $M\in\RR_{>0}^{\NN}$ satisfying

$M\le N$,

$\liminf_{p\rightarrow+\infty}(M_p/p!)^{1/p}>0$ resp. $\lim_{p\rightarrow+\infty}(M_p/p!)^{1/p}=+\infty$ and
$$\mathcal{B}(\mathcal{D}_{\{N\}}([-1,1]))\supseteq\Lambda_{\{M\}},\;\;\;\text{resp.}\;\;\;\mathcal{B}(\mathcal{D}_{(N)}([-1,1]))\supseteq\Lambda_{(M)},$$
has a maximal element (which is given by \eqref{optimalSVsequence}).
\end{theorem}

The knowledge of this optimal sequence will provide more information on the difference between the relevant conditions and give an answer to the following question: Given a nonquasianalytic sequence $N$, which (strictly smaller) sequence $M$ is maximal among all sequences allowing a mixed setting? This also gives a first piece of information related to the following problem: How far is the Borel map from being surjective in the nonquasianalytic but not strongly nonquasianalytic setting? Unfortunately, it seems that for the proofs from \cite{Borelmapgenericity} and \cite{Borelmapalgebraity} the quasianalyticity of the weight sequence is indispensable; but can we transfer the results to this situation? This question has been asked by Prof. Javier Sanz after a talk of the author about the results from \cite{Borelmapalgebraity} during a research stay at the Universidad de Valladolid.\vspace{6pt}

However, the optimal sequence for $(M,N)_{SV}$ turns out to be quite technical and involved, and so we are also interested in comparing it with the descendant which has much better regularity and growth properties. We prove that when, roughly speaking, the function $p\mapsto\frac{\nu_p}{p}\sum_{k\ge p}\frac{1}{\nu_k}$ is not ''too irregular'', then no difference occurs on the level of weighted sequence spaces between both optimal sequences, see Theorem \ref{optimalSV2}.\vspace{6pt}

For completeness let us also mention the following: In the mixed weight function setting there exists only one relevant mixed condition for the Borel map for both types, see \cite{BonetMeiseTaylorSurjectivity} (for the more general Whitney jet mapping for the Roumieu-type see \cite{whitneyextensionmixedweightfunction}, \cite{whitneyextensionmixedweightfunctionII}). However, the approach from \cite{BonetMeiseTaylorSurjectivity} has been inspired by \cite{Carleson61} and \cite{ehrenpreisbook} where the mixed weight sequence case has been treated by involving the so-called {\itshape associated weight functions}. More comments are given in Remark \ref{Carlesonremark} which has been motivated by the observations by one of the anonymous referees.

Let us also emphasize that in the literature there exist different approaches and techniques in order to study the image of the Borel map, e.g. by using summation methods. In this context we refer to the very recent paper \cite{Kiro20}; the crucial function $L$ there is connected to $M$ via $L(p)=(M_p/p!)^{1/p}$.\vspace{6pt}

We summarize now the structure of this article: After collecting and recalling necessary basic notation in Section \ref{notation}, in Section \ref{mainresults} we compute and study the optimal sequence for $(M,N)_{SV}$ (in Sect. \ref{mixedsv}) and compare it with the descendant, see Theorem \ref{optimalSV2} in Sect. \ref{descendentcompare}. In Section \ref{ramisection} we give some comments on the results of Section \ref{mainresults} for the $r$-ramified setting. In Section \ref{counterex} we construct an (involved) example which shows that in general there is a difference between the optimal sequences.\vspace{6pt}

In the final Section \ref{failureinjandsurj} we provide some information on the extent of the failure of the injectivity and surjectivity of $\mathcal{B}$ for sequences which are nonquasianalytic but not strongly nonquasianalytic, see Propositions \ref{noninj} and \ref{nqfailure}. The problem of injectivity, treated in Section \ref{failureinj}, gives information on the size of the kernel of $\mathcal{B}$ (resp. on the quantity of ultradifferentiable flat functions). This question is somehow ''dual'' to the problem of the failure of surjectivity. It has been motivated by a question asked by Prof. Fernando Sanz (also from the Universidad de Valladolid).

\vspace{6pt}

\textbf{Acknowledgements.}
The author wishes to thank Armin Rainer and David N. Nenning, both from the University of Vienna, for interesting and helpful discussions during the preparation of this article and for their careful reading of a preliminary version.\vspace{6pt}

The author also thanks the two anonymous referees for their careful reading and valuable suggestions which have improved and clarified the presentation.

\section{Notation}\label{notation}
We start by collecting all conditions on weights and ultradifferentiable spaces needed in the main results below. Throughout this article we will write $\NN:=\{0,1,2,\dots\}$ and $\NN_{>0}:=\{1,2,\dots\}$.

\subsection{Weight sequences}\label{weightsequences}
Given a sequence $M=(M_j)_j\in\RR_{>0}^{\NN}$ we also use $m=(m_j)_j$ defined by $m_j:=\frac{M_j}{j!}$ and $\mu_j:=\frac{M_j}{M_{j-1}}$, $\mu_0:=1$. The sequence $M$ is called {\itshape normalized} if $1=M_0\le M_1$, which can always be assumed without loss of generality. For any $r>0$ we denote the $r$-th power of $M$ by $M^r=(M^r_j)_{j\in\NN}$.\vspace{6pt}

$M$ is called {\itshape log-convex} if
$$\forall\;j\in\NN_{>0}:\;M_j^2\le M_{j-1} M_{j+1},$$
equivalently if $(\mu_j)_j$ is nondecreasing. If $M$ is log-convex and normalized, then both $j\mapsto M_j$ and $j\mapsto(M_j)^{1/j}$ are nondecreasing and $(M_j)^{1/j}\le\mu_j$ for all $j\in\NN_{>0}$.

If $m$ is log-convex, then $M$ is called {\itshape strongly log-convex}, denoted by \hypertarget{slc}{$(\text{slc})$}. For our purposes it is convenient to consider the set of sequences
$$\hypertarget{LCset}{\mathcal{LC}}:=\{M\in\RR_{>0}^{\NN}:\;M\;\text{is normalized, log-convex},\;\lim_{j\rightarrow+\infty}(M_j)^{1/j}=+\infty\}.$$
We see that $M\in\hyperlink{LCset}{\mathcal{LC}}$ if and only if $1=\mu_0\le\mu_1\le\dots$, $\lim_{j\rightarrow+\infty}\mu_j=+\infty$ (see e.g. \cite[p. 104]{compositionpaper}) and there is a one-to-one correspondence between $M$ and $\mu=(\mu_j)_j$ by taking $M_j:=\prod_{i=0}^j\mu_i$.\vspace{6pt}

$M$ has {\itshape moderate growth}, denoted by \hypertarget{mg}{$(\text{mg})$}, if
$$\exists\;C\ge 1\;\forall\;j,k\in\NN:\;M_{j+k}\le C^{j+k} M_j M_k.$$
In \cite{Komatsu73} it is denoted by $(M.2)$ and called {\itshape stability under ultradifferential operators.} We can replace in this condition $M$ by $m$ by changing the constants. It is known (see e.g. \cite[Lemma 2.2]{whitneyextensionweightmatrix}) that for any given $M\in\hyperlink{LCset}{\mathcal{LC}}$ condition \hyperlink{mg}{$(\on{mg})$} is equivalent to $\sup_{j\in\NN}\frac{\mu_{2j}}{\mu_j}<+\infty$ and to $\sup_{j\in\NN_{>0}}\frac{\mu_j}{(M_j)^{1/j}}<+\infty$. The latter condition shows that sequences of quotients and roots are comparable up to a constant.\vspace{6pt}

$M$ satisfies \hypertarget{gamma1}{$(\gamma_1)$} (see \cite{petzsche}) if
$$\sup_{j\in\NN_{>0}}\frac{\mu_j}{j}\sum_{k\ge j}\frac{1}{\mu_k}<+\infty.$$
In the literature \hyperlink{gamma1}{$(\gamma_1)$} is also called ''strong nonquasianalyticity'', see \cite{petzsche}, and in \cite{Komatsu73} it is denoted by $(M.3)$. In \cite{petzsche} the surjectivity of $\mathcal{B}$ has been characterized in terms of \hyperlink{gamma1}{$(\gamma_1)$}.\vspace{6pt}




Let $M,N\in\RR_{>0}^{\NN}$ be given, we write $M\hypertarget{preceq}{\preceq}N$ if $\sup_{j\in\NN_{>0}}\left(\frac{M_j}{N_j}\right)^{1/j}<+\infty$. We call $M$ and $N$ {\itshape equivalent}, written $M\hypertarget{approx}{\approx}N$, if $M\hyperlink{preceq}{\preceq}N$ and $N\hyperlink{preceq}{\preceq}M$. This equivalence preserves \hyperlink{mg}{$(\text{mg})$}. Finally, write $M\le N$ if $M_j\le N_j$ for all $j\in\NN$.\vspace{6pt}

We mention that in \cite[Prop. 1.1]{petzsche} it has been shown that \hyperlink{gamma1}{$(\gamma_1)$} for log-convex $M$ implies that there exists an equivalent sequence $N$ having \hyperlink{slc}{$(\text{slc})$}, so \hyperlink{gamma1}{$(\gamma_1)$} ''implies'' \hyperlink{slc}{$(\text{slc})$}.\vspace{6pt}

We denote by $M^{\on{lc}}$ the log-convex minorant of $M$, i.e. $M^{\on{lc}}$ is the largest sequence among all sequences $L$ which are log-convex and satisfy $L\le M$. For concrete formulas computing $M^{\on{lc}}$ we refer to \cite{hoermander}, alternatively $M^{\on{lc}}$ can be obtained by using the so-called associated weight function of $M$, see \cite[$(3.2)$, Prop. 3.2]{Komatsu73} and \cite[Chapitre I, 1.8]{mandelbrojtbook}.


\subsection{Ultradifferentiable sequence and function spaces}
Let $M\in\RR_{>0}^{\NN}$ and $h>0$ be given, for a sequence $\mathbf{a}:=(a_p)_p\in\CC^{\NN}$ we put
$$|\mathbf{a}|_{M,h}:=\sup_{p\in\NN}\frac{|a_p|}{h^p M_p},\hspace{20pt}\Lambda_{M,h}:=\{(a_p)_p\in\CC^{\NN}: |\mathbf{a}|_{M,h}<+\infty\}.$$
Furthermore we set
$$\Lambda_{(M)}:=\{(a_p)_p\in\CC^{\NN}: \forall\;h>0: |\mathbf{a}|_{M,h}<+\infty\},$$
and
$$\Lambda_{\{M\}}:=\{(a_p)_p\in\CC^{\NN}: \exists\;h>0: |\mathbf{a}|_{M,h}<+\infty\}.$$
$\Lambda_{(M)}$ is a weighted sequence space of {\itshape Beurling-type}, and $\Lambda_{\{M\}}$ of {\itshape Roumieu-type.} We endow $\Lambda_{(M)}$ resp. $\Lambda_{\{M\}}$ with a natural projective, respectively inductive, topology via $$\Lambda_{(M)}=\underset{h>0}{\varprojlim}\;\Lambda_{M,h},\hspace{30pt}\Lambda_{\{M\}}=\underset{h>0}{\varinjlim}\;\Lambda_{M,h}.$$
If $M\hyperlink{preceq}{\preceq} M'$, then clearly $\Lambda_{[M]}\subseteq\Lambda_{[M']}$, so equivalence of weight sequences preserves the associated weighted sequence spaces.\vspace{6pt}

Analogously, the (local) ultradifferentiable class of Roumieu-type is given by
\begin{align*}
&\mathcal{E}_{\{M\}}(\RR,\CC):=
\\&
\{f\in\mathcal{E}(\RR,\CC):\;\;\;\forall\;K\subseteq\RR\;\text{compact}\;\exists\;C,h>0\;\forall\;j\in\NN\;\forall\;x\in K:\;\;\;|f^{(j)}(x)|\le Ch^jM_j\},
\end{align*}
and the Beurling-type by
\begin{align*}
&\mathcal{E}_{(M)}(\RR,\CC):=
\\&
\{f\in\mathcal{E}(\RR,\CC):\;\;\;\forall\;K\subseteq\RR\;\text{compact}\;\forall\;h>0\;\exists\;C_h>0\;\forall\;j\in\NN\;\forall\;x\in K:\;\;\;|f^{(j)}(x)|\le C_hh^jM_j\}.
\end{align*}
If $M\hyperlink{preceq}{\preceq} M'$, then clearly $\mathcal{E}_{[M]}\subseteq\mathcal{E}_{[M']}$ is valid. Moreover, $\liminf_{p\rightarrow+\infty}(m_p)^{1/p}>0$ implies $\mathcal{E}_{\{M^{\on{lc}}\}}=\mathcal{E}_{\{M\}}$ and $\lim_{p\rightarrow+\infty}(m_p)^{1/p}=+\infty$ implies $\mathcal{E}_{(M^{\on{lc}})}=\mathcal{E}_{(M)}$, see \cite[Theorem 2.15]{compositionpaper}. However, in general one has only $\Lambda_{[M^{\on{lc}}]}\subseteq\Lambda_{[M]}$.\vspace{6pt}

Finally, for each $h>0$ and we define the Banach space
\begin{align*}
&\mathcal{D}_{M,h}([-1,1]):=\left\{f\in\mathcal{E}(\RR,\CC):\;\supp(f)\subseteq[-1,1], \sup_{p\in\NN, x\in\RR}\frac{|f^{(p)}(x)|}{h^pM_p}<+\infty\right\},
\end{align*}
and the ultradifferentiable test function class of Roumieu-type
\begin{equation*}\label{Roumieu-LB-space}
\mathcal{D}_{\{M\}}([-1,1]):=\underset{h>0}{\varinjlim}\;\mathcal{D}_{M,h}([-1,1]),
\end{equation*}
which is a countable $(LB)$-space, and of Beurling-type
\begin{equation*}\label{Beurling-Frechet-space}
\mathcal{D}_{(M)}([-1,1]):=\underset{h>0}{\varprojlim}\;\mathcal{D}_{M,h}([-1,1]),
\end{equation*}
which is a Fr\'{e}chet space.

Then for the {\itshape Borel map} $\mathcal{B}$ (at $0$) we get
$$\mathcal{B}:\mathcal{D}_{[M]}([-1,1]),\mathcal{E}_{[M]}(\RR,\CC)\longrightarrow\Lambda_{[M]},\hspace{20pt}\mathcal{B}(f):=(f^{(p)}(0))_{p\in\NN}.$$

The nontriviality of the classes $\mathcal{D}_{[M]}([-1,1])$ is characterized in terms of $M$ by the nonquasianalyticity condition, see \cite[Theorem 1.3.8]{hoermander} and \cite[Theorem 4.2]{Komatsu73}.

If $N\in\hyperlink{LCset}{\mathcal{LC}}$, then $\mathcal{D}_{[N]}([-1,1])\neq\{0\}$ if and only if $N$ is {\itshape nonquasianalytic,} denoted by \hypertarget{mnq}{$(\text{nq})$}, which means that
$$\sum_{p=1}^{+\infty}\frac{1}{\nu_p}<+\infty.$$
Clearly, each $N$ satisfying \hyperlink{gamma1}{$(\gamma_1)$} is nonquasianalytic and it is known that \hyperlink{mnq}{$(\text{nq})$} implies $\lim_{p\rightarrow+\infty}(n_p)^{1/p}=\lim_{p\rightarrow+\infty}\nu_p/p=+\infty$, see e.g. \cite[Lemma 4.1]{Komatsu73}, which does imply that both types of ultradifferentiable functions contain the class of real-analytic functions, see e.g. \cite[Prop. 2.12 $(4),(5)$]{compositionpaper}. Equivalence of sequences belonging to \hyperlink{LCset}{$\mathcal{LC}$} preserves \hyperlink{mnq}{$(\text{nq})$}.\vspace{6pt}

For given $N\in\hyperlink{LCset}{\mathcal{LC}}$ we introduce the sets of sequences
\begin{equation*}
\mathcal{N}_{\preceq,R}:=\{M\in\RR_{>0}^{\NN}:\;\;\;\liminf_{p\rightarrow+\infty}(m_p)^{1/p}>0,\;\;\;M\hyperlink{preceq}{\preceq}N\},
\end{equation*}
\begin{equation*}
\mathcal{N}_{\preceq,B}:=\{M\in\RR_{>0}^{\NN}:\;\;\;\lim_{p\rightarrow+\infty}(m_p)^{1/p}=+\infty,\;\;\;M\hyperlink{preceq}{\preceq}N\}.
\end{equation*}
Obviously $\mathcal{N}_{\preceq,B}\subseteq\mathcal{N}_{\preceq,R}$; $N\in\mathcal{N}_{\preceq,R}$ provided that $\liminf_{p\rightarrow+\infty}(n_p)^{1/p}>0$; and finally $N\in\mathcal{N}_{\preceq,B}$ provided that $\lim_{p\rightarrow+\infty}(n_p)^{1/p}=+\infty$. We also introduce the smaller sets
\begin{equation*}
\mathcal{N}_{\preceq,\mathcal{LC},R}:=\{M\in\mathcal{N}_{\preceq,R}:\;\;\;M\in\hyperlink{LCset}{\mathcal{LC}}\},
\end{equation*}
and
\begin{equation*}
\mathcal{N}_{\preceq,\mathcal{LC},B}:=\{M\in\mathcal{N}_{\preceq,B}:\;\;\;M\in\hyperlink{LCset}{\mathcal{LC}}\}.
\end{equation*}
Note that by normalization the relation $M\hyperlink{preceq}{\preceq}N$ means precisely $M_p\le C^pN_p$ for some $C\ge 1$ and all $p\in\NN$. Hence replacing $M$ by the equivalent sequence $\widetilde{M}^C:=(M_p/C^p)_{p\in\NN}$ we get $\widetilde{M}^C\le N$. This shows that, since both the weighted sequence spaces $\Lambda_{[M]}$ and the ultradifferentiable (test) function classes are stable w.r.t. \hyperlink{approx}{$\approx$}, for our purposes instead of treating $\mathcal{N}_{\preceq,R}$, $\mathcal{N}_{\preceq,B}$, equivalently we could also consider the set
\begin{equation*}\label{MbelowNroum}
\mathcal{N}_{\le,R}:=\{M\in\RR_{>0}^{\NN}:\;\;\;\liminf_{p\rightarrow+\infty}(m_p)^{1/p}>0,\;\exists\;C\ge 1:\;\;\;M\le CN\},
\end{equation*}
resp.
\begin{equation*}\label{MbelowNbeur}
\mathcal{N}_{\le,B}:=\{M\in\RR_{>0}:\;\;\;\lim_{p\rightarrow+\infty}(m_p)^{1/p}=+\infty,\;\exists\;C\ge 1:\;\;\;M\le CN\},
\end{equation*}
and similarly $\mathcal{N}_{\le,\mathcal{LC},R}$, $\mathcal{N}_{\le,\mathcal{LC},B}$, see also \cite[Remark 3.1 $(i)$]{mixedramisurj}. Analogously also the relevant mixed conditions are not affected by such a modification, see $(iii)$ in Remark \ref{lemma24remark}. Note that $M$ is (strongly) log-convex if and only if $\widetilde{M}^C$ is so for some/each $C>0$.

\subsection{Strong nonquasianalyticity conditions in the mixed setting}\label{mixedsnqs}
For any given nonquasianalytic $N\in\hyperlink{LCset}{\mathcal{LC}}$ and $M\in\RR_{>0}^{\NN}$ we now recall the mixed conditions: First $(M,N)_{SV}$ means
$$\exists\;s\in\NN_{>0}:\;\;\sup_{p\in\NN_{>0}}\frac{\lambda_{p,s}^{M,N}}{p}\sum_{k\ge p}\frac{1}{\nu_k}<+\infty,$$
with $\lambda^{M,N}_{p,s}:=\sup_{0\le j<p}\left(\frac{M_p}{s^p N_j}\right)^{1/(p-j)}$, and second
$$\sup_{p\in\NN_{>0}}\frac{\mu_p}{p}\sum_{k\ge p}\frac{1}{\nu_k}<+\infty,$$
denoted by $(M,N)_{\gamma_1}$. Next we summarize consequences for these conditions when $M$ belongs to some sets introduced before.

\begin{remark}\label{lemma24remark}
Inspecting the proof of \cite[Lemma 2.4]{mixedramisurj} (see also \cite[$2(a)$]{surjectivity}) we get the following (translated into the notation from the previous section):
\begin{itemize}
\item[$(i)$] If $M\in\mathcal{N}_{\le,R}$ resp. $M\in\mathcal{N}_{\le,B}$, then
\begin{equation}\label{lemma24remarkequ}
\forall\;p,s\in\NN_{>0}:\;\;\;\lambda^{M,N}_{p,s}\le C\nu_p,
\end{equation}
$C$ denoting the constant from $M\le CN$.

\item[$(ii)$] If $M\in\mathcal{N}_{\le,\mathcal{LC},R}$ resp. $M\in\mathcal{N}_{\le,\mathcal{LC},B}$, then even $$\forall\;p,s\in\NN_{>0}:\;\;\;\lambda^{M,N}_{p,s}\le C\min\{\mu_p,\nu_p\},$$
    with $C$ denoting the constant from $M\le CN$. Hence in this situation $(M,N)_{\gamma_1}$ implies $(M,N)_{SV}$. If $M$ in addition satisfies \hyperlink{mg}{$(\on{mg})$}, which is in this case equivalent to $\sup_{j\in\NN_{>0}}\frac{\mu_j}{(M_j)^{1/j}}<+\infty$, then also the converse is true and so $(M,N)_{\gamma_1}$ is equivalent to $(M,N)_{SV}$.

\item[$(iii)$] $(M,N)_{\gamma_1}$ resp. $(M,N)_{SV}$ is valid iff $(\widetilde{M}^C,N)_{\gamma_1}$ resp. $(\widetilde{M}^C,N)_{SV}$ holds for some/any $\widetilde{M}^C:=(M_p/C^p)_{p\in\NN}$, $C\ge 1$ (we can take $C\in\NN_{>0}$). The parameter $s\in\NN_{>0}$ in $(M,N)_{SV}$ is then multiplied by $C$.
\end{itemize}
\end{remark}

Now we recall the main results \cite[Thm. 3.2, Thm. 4.2]{mixedramisurj} (with $r=1$ there) in our notation, see also \cite{surjectivity}.

\begin{theorem}\label{SVtheorem}
Let $N\in\hyperlink{LCset}{\mathcal{LC}}$ and $M\in\mathcal{N}_{\preceq,R}$ resp. $M\in\mathcal{N}_{\preceq,B}$. Then the following properties are equivalent:
\begin{itemize}
\item[$(i)$] $\mathcal{B}(\mathcal{D}_{\{N\}}([-1,1]))\supseteq\Lambda_{\{M\}}$ resp. $\mathcal{B}(\mathcal{D}_{(N)}([-1,1]))\supseteq\Lambda_{(M)}$,

\item[$(ii)$] $(M,N)_{SV}$ is valid.
\end{itemize}
Consequently, if $(i)$ or $(ii)$ holds true, then $N$ has to be nonquasianalytic.
\end{theorem}

In fact, in \cite{mixedramisurj} the assumption on $M$ has been $M\in\mathcal{N}_{\preceq,\mathcal{LC},R}$ resp. $M\in\mathcal{N}_{\preceq,\mathcal{LC},B}$ and we point out:

\begin{itemize}
\item[$(a)$] For the Roumieu case, a careful inspection of the proof of \cite[Thm. 3.2]{mixedramisurj} (with $r=1$) shows that the assumption log-convexity on $M$ is superfluous and the estimate \eqref{lemma24remarkequ} suffices to conclude.

\item[$(b)$] For the Beurling case \cite[Thm. 4.2]{mixedramisurj} we give more details: In the proof of \cite[Thm. 4.5]{mixedramisurj} the technical result of \cite[Lemme 16]{ChaumatChollet94} has been applied in order to reduce this situation to the Roumieu case. We have to avoid the choice $\gamma_k:=\frac{1}{\mu_k}$ since $k\mapsto\mu_k$ is not necessarily nondecreasing any more. In order to conclude, take e.g. $\gamma_k:=\frac{1}{k}$ or any other nonincreasing sequence of positive real numbers tending to $0$ as $k\rightarrow+\infty$. Then the constructed sequence $R$ will not necessarily be log-convex, i.e. $(i)$ on \cite[p. 562]{mixedramisurj} will fail. However, the other properties of the sequences $R$ and $S$ listed there are still valid and sufficient to conclude since, as pointed out above for the Roumieu case, the log-convexity of the smaller sequence is not required necessarily.

\item[$(c)$] The same is true for the more general ramified case mentioned in in Section \ref{ramisection}.
\end{itemize}

We close this section with the following observation. It has been motivated by a question asked by one of the anonymous referees who has also pointed out the literature references \cite{Carleson61} and \cite{ehrenpreisbook} and their strong connection to the topics and problems studied in this article.

\begin{remark}\label{Carlesonremark}
Indeed, we have a third mixed relevant condition for weight sequences; not given directly by $M$ and $N$ but by involving the so-called {\itshape associated weight functions.}
\begin{itemize}
\item[$(i)$] Given $M\in\hyperlink{LCset}{\mathcal{LC}}$, the {\itshape associated function} $\omega_M: \RR_{\ge 0}\rightarrow\RR\cup\{+\infty\}$ is defined by
\begin{equation*}\label{assofunc}
\omega_M(t):=\sup_{p\in\NN}\log\left(\frac{t^p}{M_p}\right)\;\;\;\text{for}\;t>0,\hspace{30pt}\omega_M(0):=0,
\end{equation*}
we refer to \cite[Chapitre I]{mandelbrojtbook}, see also \cite[Definition 3.1]{Komatsu73}. The given sequence can be expressed in terms of $\omega_M$ by
\begin{equation}\label{Mvsassofunc}
M_p=\sup_{t\ge 0}\frac{t^p}{\exp(\omega_{M}(t))},\;\;\;p\in\NN.
\end{equation}
With this notation the mixed condition reads
\begin{equation}\label{gammarfctmix}
(\omega_M,\omega_N)_{\on{snq}}:\Leftrightarrow\exists C>0\;\forall t\ge 0 :\int_1^{+\infty}\frac{\omega_N(ty)}{y^{2}}dy\le C\omega_M(t)+C.
\end{equation}
\item[$(ii)$] \eqref{gammarfctmix} has been used in \cite{Carleson61} for the Roumieu-type and in \cite[Sect. XIII.3]{ehrenpreisbook} for the Beurling-type in order to treat extensions $(i)$ in Theorem \ref{SVtheorem}. And later, when taking general weight functions (in the sense of Braun-Meise-Taylor) $\sigma\equiv\omega_M$ and $\omega\equiv\omega_N$, for this setting in \cite{BonetMeiseTaylorSurjectivity} a precise characterization in the spirit of Theorem \ref{SVtheorem} has been given.

\item[$(iii)$] In \cite{Carleson61} only the sufficiency of \eqref{gammarfctmix} is obtained, in \cite[Theorem 13.17]{ehrenpreisbook} also the necessity is shown, and the proofs in \cite{Carleson61} and \cite[Theorem 13.17]{ehrenpreisbook} are using different methods, see also the discussion in \cite[Rem. 13.8, p. 487]{ehrenpreisbook}.

However, for \cite[Theorem 13.17]{ehrenpreisbook} in the preparatory result \cite[Lemma 13.16]{ehrenpreisbook} on p. 480 an assumption is made on the associated weight function. In fact, by \cite[Prop. 3.6]{Komatsu73} this requirement precisely means \hyperlink{mg}{$(\on{mg})$} for $N$ which is not needed in our general considerations in Sect. \ref{mixedsv} below. Therefore note that the function $\lambda$ in \cite{ehrenpreisbook} is corresponding to $\exp(\omega_N)$ in our notation, see \cite[p. 163]{ehrenpreisbook} and \cite[Thm. 19.11]{rudin} (and take also into account the equivalence shown in \cite[Lemma 3.3]{BonetMeiseTaylorSurjectivity}).

Requirement \hyperlink{mg}{$(\on{mg})$} seems to be indispensable when involving weight function techniques in order to prove results for the weight sequence setting. This fact is indicated by the comparison results obtained in \cite{BonetMeiseMelikhov07}, see also \cite{compositionpaper}.

\item[$(iv)$] In any case, it seems to be a natural and interesting question how \eqref{gammarfctmix} is related to conditions $(M,N)_{\gamma_1}$, $(M,N)_{SV}$, how the optimal (minimal) weight $\kappa_{\omega_N}(t):=\int_1^{+\infty}\frac{\omega_N(ty)}{y^{2}}dy$ (see again \cite{BonetMeiseTaylorSurjectivity}) can be used to compute a sequence $K$ by involving \eqref{Mvsassofunc} and how $K$ is in general connected to the optimal sequences for $(M,N)_{\gamma_1}$, $(M,N)_{SV}$.

    Note that in \cite[Lemma 5.7]{whitneyextensionmixedweightfunction} it is shown that for $M,N\in\hyperlink{LCset}{\mathcal{LC}}$ with $\mu\le\nu$ we have that $(M,N)_{\gamma_1}$ implies $(\omega_M,\omega_N)_{\on{snq}}$, see also \cite[Lemma 4, Rem. 6]{mixedsectorialextensions} (for $r=1$).\vspace{6pt}

    Motivated by these observations the author has started to study this question by dealing with abstractly given weight matrices $\mathcal{M}$; this will be considered in another paper.
\end{itemize}
\end{remark}

\section{Optimal sequences in the mixed setting}\label{mainresults}
\subsection{The mixed strong nonquasianalyticity condition}\label{descendant}
In \cite[Section 4.1]{whitneyextensionweightmatrix}, which is based on an idea arising in the proof of \cite[Proposition 1.1]{petzsche}, it has been shown that to each $N\in\hyperlink{LCset}{\mathcal{LC}}$ satisfying \hyperlink{mnq}{$(\on{nq})$} we can associate a sequence $S^N$ with good regularity properties denoted by the {\itshape descendant}.\vspace{6pt}

We recall the construction of $S^N$ defined by its quotients
$$\sigma^N_p:=\frac{\tau_1 p}{\tau_p},\;\;\;p\in\NN_{>0},\hspace{20pt}\sigma^N_0:=1,$$
with
$$\tau_p:=\frac{p}{\nu_p}+\sum_{j\ge p}\frac{1}{\nu_j},\;\;\;p\ge 1,$$
and $S^N$ satisfies the following properties (see \cite[Lemma 4.2]{whitneyextensionweightmatrix}):

\begin{itemize}
\item[$(i)$] $\sigma^N_p\ge 1$ for all $p\in\NN$ and $s^N:=(S^N_p/p!)_{p\in\NN}\in\hyperlink{LCset}{\mathcal{LC}}$ (so $S^N$ is strongly log-convex),

\item[$(ii)$] there exists $C>0$ such that $\sigma^N_p\le C\nu_p$ for all $p\in\NN$,

\item[$(iii)$] $(S^N,N)_{\gamma_1}$,

\item[$(iv)$] if $N$ enjoys \hyperlink{mg}{$(\on{mg})$}, then so does $S^N$,

\item[$(v)$] $S^N$ is optimal/maximal in the following sense: If $M\in\hyperlink{LCset}{\mathcal{LC}}$ is given with $\mu_p\le C\nu_p$ for some $C\ge 1$ and $(M,N)_{\gamma_1}$, then $\mu_p\le D\sigma^N_p$ follows for some constant $D\ge 1$.

We also have that
$$\exists\;B\ge 1\;\forall\;p\in\NN:\;\;\;B^{-1}\sigma^N_p\le\nu_p\le B\sigma_p^N,$$
(which implies $S^N\hyperlink{approx}{\approx}N$) if and only if $N$ satisfies \hyperlink{gamma1}{$(\gamma_1)$} resp. if and only if $(N,N)_{SV}$, see \cite[Theorem 5.2]{mixedramisurj}.
\end{itemize}

\begin{remark}\label{desremark}
By the previous comments $(i)$ and $(ii)$, $S^N\in\mathcal{N}_{\preceq,\mathcal{LC},B}$. We can replace $S^N$ by the equivalent sequence (''modified descendant'') $\widetilde{S^N}$ defined by its quotients as follows: We put $\widetilde{\sigma^N}_p:=\sigma^N_p/C$ for $C$ denoting the constant from $(ii)$ (we can assume $C\in\NN_{>0}$) and for all $p\ge p_C$, with $p_C\in\NN$ chosen minimal to ensure $\sigma^N_p\ge C$. This is possible by having $\lim_{p\rightarrow+\infty}\sigma^N_p=+\infty$. Finally, for all $1\le p<p_C$ we put $\widetilde{\sigma^N}_p:=1(=\sigma^N_1\le\nu_1)$.

So $\widetilde{\sigma^N}\le\nu$, $\widetilde{S^N}\in\mathcal{N}_{\le,\mathcal{LC},B}$ and $(\widetilde{S^N},N)_{\gamma_1}$ holds by combining comment $(iii)$ above and $(iii)$ in Remark \ref{lemma24remark}. By $(ii)$ in Remark \ref{lemma24remark} also $(\widetilde{S^N},N)_{SV}$ holds true and by $(iii)$ there finally property $(S^N,N)_{SV}$ follows.


This shows that for any $N\in\hyperlink{LCset}{\mathcal{LC}}$ satisfying \hyperlink{mnq}{$(\on{nq})$} the set of all sequences belonging to $\mathcal{N}_{\preceq,\mathcal{LC},B}$ (hence to $\mathcal{N}_{\preceq,\mathcal{LC},R}$) and satisfying $(\cdot,N)_{SV}$ is never empty.

(Of course, $N\in\mathcal{N}_{\preceq,\mathcal{LC},B}$ is valid automatically, but $(N,N)_{SV}$ if and only if $N$ satisfies \hyperlink{gamma1}{$(\gamma_1)$}.)
\end{remark}

\subsection{Optimal sequence for the mixed Schmets-Valdivia-condition}\label{mixedsv}
The aim is now to determine the maximal sequence $L$ belonging to the set $\mathcal{N}_{\preceq,R}$ resp. $\mathcal{N}_{\preceq,B}$ and satisfying $(L,N)_{SV}$, i.e. the maximal sequence such that Theorem \ref{SVtheorem} is valid and hence admitting the maximal possible control of loss of regularity within the ultradifferentiable weight sequence setting. So $L$ will determine the largest possible sequence space $\Lambda_{[L]}$ which is contained in the image of the Borel map $\mathcal{B}$ when being considered on $\mathcal{D}_{[N]}([-1,1])$ (with $N$ fixed).

\begin{lemma}\label{optimalSV}
Let $N\in\hyperlink{LCset}{\mathcal{LC}}$ and assume that $N$ is {\itshape nonquasianalytic}. For any $s\in\NN_{>0}$ define the sequence $L^s=(L^s_p)_{p\in\NN}$ as follows:
\begin{equation}\label{optimalSVsequence}
L^s_0:=1,\hspace{30pt}L^s_p:=s^p\min_{0\le j\le p-1}\left\{\left(\frac{p}{\sum_{k\ge p}\frac{1}{\nu_k}}\right)^{p-j}N_j\right\},\;\;p\in\NN_{>0}.
\end{equation}
For convenience we will always write $L\equiv L^1$.

Then:
\begin{itemize}
\item[$(i)$] $(L^s,N)_{SV}$ holds true for all $s\in\NN_{>0}$ (note that $L^s\hyperlink{approx}{\approx}L^t$ for all $s,t\in\NN_{>0}$),

\item[$(ii)$] $M\hyperlink{preceq}{\preceq}L^s$ for any $M\in\RR_{>0}^{\NN}$ satisfying $(M,N)_{SV}$ and for all $s\in\NN_{>0}$, 

\item[$(iii)$] $L^s\hyperlink{preceq}{\preceq}N$ for all $s\in\NN_{>0}$.

\end{itemize}

Since there exists $M\in\mathcal{N}_{\preceq,\mathcal{LC},R}$ resp. $M\in\mathcal{N}_{\preceq,\mathcal{LC},B}$ with $(M,N)_{SV}$, e.g. take the descendant $S^N$ (see Remark \ref{desremark}), assertions $(i)$, $(ii)$ and $(iii)$ hold true for the log-convex minorant $(L^s)^{\on{lc}}$ as well (with $(ii)$ restricted to log-convex sequences $M$).

\end{lemma}

\demo{Proof}
$(i)$ Let $s,p\in\NN_{>0}$ and $0\le l\le p-1$. Then
$$\left(\frac{L^s_p}{s^p N_l}\right)^{1/(p-l)}=\left(\frac{s^p\min_{0\le j\le p-1}\left\{\left(\frac{p}{\sum_{k\ge p}\frac{1}{\nu_k}}\right)^{p-j}N_j\right\}}{s^p N_l}\right)^{1/(p-l)}\underbrace{\le}_{l=j}\frac{p}{\sum_{k\ge p}\frac{1}{\nu_k}},$$
hence $(L^s,N)_{SV}$ follows.\vspace{6pt}

$(ii)$ $(M,N)_{SV}$ precisely means that
$$\exists\;s_0\in\NN_{>0}\;\exists\;C\ge 1\;\forall\;p\in\NN_{>0}\;\forall\;0\le j\le p-1:\;\;\;M_p\le s_0^pC^{p-j}N_j\left(\frac{p}{\sum_{k\ge p}\frac{1}{\nu_k}}\right)^{p-j}.$$
So $M_p\le C^pL^{s_0}_p$ follows for all $p\in\NN_{>0}$ for this value $s_0\in\NN_{>0}$. Since $L^s\hyperlink{approx}{\approx}L^t$ for all $s,t\in\NN_{>0}$, we are done.\vspace{6pt}

$(iii)$ By choosing $j=p-1$ we have
$$\forall\;p\in\NN_{>0}:\;\;\;L^s_p\le\frac{p}{\sum_{k\ge p}\frac{1}{\nu_k}}s^pN_{p-1},$$
and so for $L^s\hyperlink{preceq}{\preceq}N$ it is sufficient to show
$$\exists\;C\ge 1\;\forall\;p\in\NN_{>0}:\;\;\;\frac{p}{\sum_{k\ge p}\frac{1}{\nu_k}}s^pN_{p-1}\le C^pN_p,$$
or equivalently
$$\exists\;C\ge 1\;\forall\;p\in\NN_{>0}:\;\;\;\left(\frac{s}{C}\right)^p\le\frac{\nu_p}{p}\sum_{k\ge p}\frac{1}{\nu_k}.$$
For $\varepsilon>0$ small enough, and for all $p\in\NN_{>0}$,
$$\frac{\nu_p}{p}\sum_{k\ge p}\frac{1}{\nu_k}=\frac{1}{p}+\frac{\nu_p}{p}\sum_{k\ge p+1}\frac{1}{\nu_k}\ge\frac{1}{p}\ge\varepsilon^p,$$
since $\frac{1}{p^{1/p}}\rightarrow 1$ as $p\rightarrow+\infty$. Hence the desired estimate is valid by choosing $C$ sufficiently large (depending on $N$ and given $s$).\vspace{6pt}

Given $M\in\mathcal{N}_{\preceq,\mathcal{LC},R}$ resp. $M\in\mathcal{N}_{\preceq,\mathcal{LC},B}$ with $(M,N)_{SV}$ by $(ii)$ we have $\liminf_{p\rightarrow+\infty}(L^s_p/p!)^{1/p}>0$ resp. $\lim_{p\rightarrow+\infty}(L^s_p/p!)^{1/p}=+\infty$ for all $s\in\NN_{>0}$. $(i)$ and $(iii)$ for $(L^s)^{\on{lc}}$ follow from $(L^s)^{\on{lc}}\le L^s$. For $(ii)$ we apply \cite[Lemma 2.6]{compositionpaper}: here $M^{\on{lc}}=M\hyperlink{preceq}{\preceq}L^s$ implies $M^{\on{lc}}=M\hyperlink{preceq}{\preceq}(L^s)^{\on{lc}}$.
\qed\enddemo

\begin{remark}\label{firstcomparisonrem}
Let $N\in\hyperlink{LCset}{\mathcal{LC}}$ be a nonquasianalytic sequence.\vspace{6pt}

By definition we get
$$\frac{\nu_p}{\sigma_p}=\frac{\nu_p/p}{\sigma_p/p}=\frac{1}{\tau_1}\frac{\nu_p}{p}\tau_p=\frac{1}{\tau_1}\left(1+\frac{\nu_p}{p}\sum_{k\ge p}\frac{1}{\nu_k}\right),$$
which means that the growth of $p\mapsto\frac{\nu_p}{p}\sum_{k\ge p}\frac{1}{\nu_k}$ measures and restricts the size of the difference between $N$ and its descendant $S^N$.

In particular, when $N$ satisfies \hyperlink{gamma1}{$(\gamma_1)$}, then $\sup_{p\in\NN_{>0}}\frac{\nu_p}{\sigma_p}<+\infty$ and so $S^N\hyperlink{approx}{\approx}N$ (as expected).\vspace{6pt}

Similarly, for the sequences $L^s$ we get when taking $j=0$ in \eqref{optimalSVsequence} that
$$\left(\frac{N_p}{L^s_p}\right)^{1/p}\ge\frac{1}{s}\frac{(N_p)^{1/p}}{p}\sum_{k\ge p}\frac{1}{\nu_k},$$
so here the growth behavior of $p\mapsto\frac{(N_p)^{1/p}}{p}\sum_{k\ge p}\frac{1}{\nu_k}$ can be used to determine the difference between $N$ and $L^s$. Note that
$$\frac{(N_p)^{1/p}}{p}\sum_{k\ge p}\frac{1}{\nu_k}\le\frac{\nu_p}{p}\sum_{k\ge p}\frac{1}{\nu_k}$$
is valid by log-convexity.
\end{remark}

Next we summarize some more consequences of Lemma \ref{optimalSV}.

\begin{itemize}
\item[$(a)$] For any given $M\in\mathcal{N}_{\preceq,R}$ resp. $M\in\mathcal{N}_{\preceq,B}$ with $(M,N)_{SV}$
  we see that $M\hyperlink{preceq}{\preceq}L^s$ for some/each $s\in\NN_{>0}$. By Remark \ref{desremark}, this holds true for the descendant $S^N\equiv M$, which implies $L^s\in\mathcal{N}_{\preceq,B}\subseteq\mathcal{N}_{\preceq,R}$, too. Consequently, some/each $L^s$ is \hyperlink{preceq}{$\preceq$}-maximal among all $M\in\mathcal{N}_{\preceq,B}$ (resp. among all $M\in\mathcal{N}_{\preceq,R}$) and having $(M,N)_{SV}$, i.e. $L^s$ is the maximal sequence such that Theorem \ref{SVtheorem} can be applied.

\item[$(b)$] The same holds true for $(L^s)^{\on{lc}}$ instead of $L^s$ by using again $S^N$ and $(L^s)^{\on{lc}}$. Note that formally $(L^s)^{\on{lc}}\in\mathcal{N}_{\preceq,\mathcal{LC},B}$ is not clear since normalization might fail for $(L^s)^{\on{lc}}$. However, $(L^s)^{\on{lc}}$ is equivalent to a sequence $\widetilde{(L^s)^{\on{lc}}}\in\mathcal{N}_{\preceq,\mathcal{LC},B}$, which is \hyperlink{preceq}{$\preceq$}-maximal among all $M\in\mathcal{N}_{\preceq,\mathcal{LC},B}$ resp. $M\in\mathcal{N}_{\preceq,\mathcal{LC},B}$.

\end{itemize}

By combining Theorem \ref{SVtheorem}, Remark \ref{desremark} and Lemma \ref{optimalSV} we immediately get the following result:

\begin{theorem}\label{maximalelement}
Let $N\in\hyperlink{LCset}{\mathcal{LC}}$. Then the set of sequences $M\in\mathcal{N}_{\preceq,R}$ resp. $M\in\mathcal{N}_{\preceq,B}$ and satisfying
$$\mathcal{B}(\mathcal{D}_{\{N\}}([-1,1]))\supseteq\Lambda_{\{M\}},\;\;\;\text{resp.}\;\;\;\mathcal{B}(\mathcal{D}_{(N)}([-1,1]))\supseteq\Lambda_{(M)},$$
has a maximal element which is given by some/each $L^s$.\vspace{6pt}

Alternatively, we can use $\widetilde{(L^s)^{\on{lc}}}$ which is \hyperlink{preceq}{$\preceq$}-maximal among all $M\in\mathcal{N}_{\preceq,\mathcal{LC},R}$ resp. $M\in\mathcal{N}_{\preceq,\mathcal{LC},B}$.
\end{theorem}

Note that for this result it is not necessary to assume that $N$ is {\itshape nonquasianalytic}; more precisely we give the following remark:

\begin{remark}\label{lquasianalytic}
Let $N\in\hyperlink{LCset}{\mathcal{LC}}$.
\begin{itemize}
\item[$(a)$] If there exists $M\in\mathcal{N}_{\preceq,R}$ resp. $M\in\mathcal{N}_{\preceq,B}$ satisfying $(M,N)_{SV}$, then nonquasianalyticity for $N$ follows from the very definition of property $(M,N)_{SV}$.

\item[$(b)$] On the contrary, if there exists $M\in\mathcal{N}_{\preceq,R}$ resp. $M\in\mathcal{N}_{\preceq,B}$ and satisfying $$\mathcal{B}(\mathcal{D}_{\{N\}}([-1,1]))\supseteq\Lambda_{\{M\}}\;\;\;\text{resp.}\;\;\;\mathcal{B}(\mathcal{D}_{(N)}([-1,1]))\supseteq\Lambda_{(M)},$$
 then $N$ has to be nonquasianalytic: Otherwise, by the {\itshape Denjoy-Carleman theorem} we would get $\mathcal{D}_{[N]}([-1,1])=\{0\}$, but $\Lambda_{[M]}$ contains at least all sequences with finitely many nonvanishing entries, a contradiction (the analogous argument has been given in \cite{mixedramisurj}).
\end{itemize}
\end{remark}

By involving a second parameter $C$ we get some more information on the technical sequence(s) $L^s$.

\begin{remark}\label{commentsonLs}
In the following let $C\ge 1$.

\begin{itemize}
\item[$(i)$] First we observe that instead of $L^s$ given by \eqref{optimalSVsequence} we can consider for any $C>1$ the sequence
\begin{equation}\label{optimalSVsequence1}
L^{s,C}_0:=1,\hspace{30pt}L^{s,C}_p:=s^p\min_{0\le j\le p-1}\left\{\left(\frac{Cp}{\sum_{k\ge p}\frac{1}{\nu_k}}\right)^{p-j}N_j\right\},\;\;\;p\in\NN_{>0}.
\end{equation}
Then $L^s_p\le L^{s,C}_p\le C^pL^s_p$, i.e. $L^s\hyperlink{approx}{\approx}L^{s,C}$, and the conclusions shown in Lemma \ref{optimalSV} are valid for $L^{s,C}$ in place of $L^s$ as well.

\item[$(ii)$] The expression appearing in the minimum in \eqref{optimalSVsequence1} is nondecreasing if and only if $$\left(\frac{Cp}{\sum_{k\ge p}\frac{1}{\nu_k}}\right)^{p-j-1}N_{j+1}\ge\left(\frac{Cp}{\sum_{k\ge p}\frac{1}{\nu_k}}\right)^{p-j}N_j,$$
    which is equivalent to $\nu_{j+1}\ge\frac{Cp}{\sum_{k\ge p}\frac{1}{\nu_k}}$.

Consequently, the minimum is attained either at $j\in\NN$ minimal such that $\nu_{j+1}>\frac{Cp}{\sum_{k\ge p}\frac{1}{\nu_k}}$ is valid, if this happens for some $j+1\le p-1\Leftrightarrow j\le p-2$, or at $j=p-1$ if $\frac{\nu_{p-1}}{p}\sum_{k\ge p}\frac{1}{\nu_k}\le C$.

\item[$(iii)$] If $N\in\hyperlink{LCset}{\mathcal{LC}}$ satisfies \hyperlink{gamma1}{$(\gamma_1)$} (i.e. $(N,N)_{\gamma_1}$ or equivalently $(N,N)_{SV}$, see \cite[Theorem 5.2]{mixedramisurj}), then by choosing $C$ sufficiently large, i.e. take
    $$C\ge\sup_{p\in\NN_{>0}}\frac{\nu_p}{p}\sum_{k\ge p}\frac{1}{\nu_k},$$
    we will note for all $p\in\NN_{>0}$ that above the second case holds true, which means that
    \begin{equation}\label{optimalSVsequence2}
    L^{s,C}_p=\frac{Cp}{\sum_{k\ge p}\frac{1}{\nu_k}}s^pN_{p-1},\;\;\;\forall\;p\in\NN_{>0}.
    \end{equation}
\end{itemize}
\end{remark}

By using this last observation we can prove the following result which shows that our approach is consistent with the characterization in \cite{petzsche}.

\begin{lemma}\label{gamma1}
Let $N\in\hyperlink{LCset}{\mathcal{LC}}$.
\begin{itemize}
\item[$(i)$] If $N$ satisfies \hyperlink{gamma1}{$(\gamma_1)$}, i.e. $N$ is strongly nonquasianalytic, then $N\hyperlink{approx}{\approx}L^{s,C}$ for all $s\in\NN_{>0}$ and $C\ge 1$.

\item[$(ii)$] Conversely, if $N\hyperlink{approx}{\approx}L^s$ for some/each $s\in\NN_{>0}$, then $N$ satisfies \hyperlink{gamma1}{$(\gamma_1)$}.
\end{itemize}
\end{lemma}

\demo{Proof}
$(i)$ By $(i)$ and $(iii)$ in Lemma \ref{optimalSV} it remains to show $N\hyperlink{preceq}{\preceq}L^{s,C}$. By \eqref{optimalSVsequence2} this means that, when $C>1$ is chosen sufficiently large, i.e.
$$C\ge\sup_{p\in\NN_{>0}}\frac{\nu_p}{p}\sum_{k\ge p}\frac{1}{\nu_k},$$
then for some $D\ge 1$ and all $p\in\NN_{>0}$ we want to have
$$D^pL^{s,C}_p=D^p\frac{Cp}{\sum_{k\ge p}\frac{1}{\nu_k}}s^pN_{p-1}\ge N_p\Leftrightarrow C(Ds)^p\ge\frac{\nu_p}{p}\sum_{k\ge p}\frac{1}{\nu_k}.$$
This last estimate is clearly satisfied for any $D\ge 1$ and $s\in\NN_{>0}$ by \hyperlink{gamma1}{$(\gamma_1)$} and the choice of $C$. Since $L^{t,C_1}\hyperlink{approx}{\approx}L^{s,C}$ for all $C,C_1\ge 1$ and $s,t\in\NN_{>0}$ we are done.\vspace{6pt}

$(ii)$ First, by the equivalence we see that $N$ has to be nonquasianalytic: Otherwise, the definition of $L^s$ given in \eqref{optimalSVsequence} would fail (one could set $L^s_p=0$ for any $p\ge 1$ which would be consistent with Remark \ref{lquasianalytic} above but which would contradict $N\in\hyperlink{LCset}{\mathcal{LC}}$). Second, this equivalence implies $\lim_{p\rightarrow+\infty}(L^s_p/p!)^{1/p}=+\infty$ because $\lim_{p\rightarrow+\infty}(n_p)^{1/p}=+\infty$ holds true by the nonquasianalyticity of $N$. Then \cite[Thm. 2.15]{compositionpaper} yields $\mathcal{E}_{[(L^s)^{\on{lc}}]}=\mathcal{E}_{[L^s]}=\mathcal{E}_{[N]}$ and this identity implies $N\hyperlink{approx}{\approx}(L^s)^{\on{lc}}$. (Alternatively this equivalence follows by using \cite[Lemma 2.6]{compositionpaper} and the assumption $N\hyperlink{approx}{\approx}L^s$.) Lemma \ref{optimalSV} yields $((L^s)^{\on{lc}},N)_{SV}$ and Theorem \ref{SVtheorem} implies $(N,N)_{SV}$ because the verified equivalence implies $\Lambda_{[N]}=\Lambda_{[(L^s)^{\on{lc}}]}$. Thus \hyperlink{gamma1}{$(\gamma_1)$} for $N$ follows, see \cite[Theorem 5.2]{mixedramisurj}.
\qed\enddemo

By Lemmas \ref{optimalSV}, \ref{gamma1}, Remark \ref{commentsonLs} and the comments on $S^N$ we find that $S^N\hyperlink{approx}{\approx}L^{s,C}\hyperlink{approx}{\approx}N$ if and only if $N$ satisfies \hyperlink{gamma1}{$(\gamma_1)$}.

The next statement shows that $L$ can be ''relatively'' near given $N$ even if \hyperlink{gamma1}{$(\gamma_1)$} is violated for this sequence.

\begin{lemma}\label{Lbig}
Let $N\in\hyperlink{LCset}{\mathcal{LC}}$ and assume that
\begin{equation}\label{Lbigequ}
\liminf_{p\rightarrow+\infty}\frac{\nu_p}{p}\sum_{k\ge p}\frac{1}{\nu_k}<1.
\end{equation}
Then $N_p\le L_p$ is valid for infinitely many numbers $p\in\NN_{>0}$.
\end{lemma}

\demo{Proof}
By $(ii)$ in Remark \ref{commentsonLs} (with $C=1$) we get $L_p=\frac{p}{\sum_{k\ge p}\frac{1}{\nu_k}}N_{p-1}$ for infinitely many values $p\in\NN_{>0}$ (see \eqref{optimalSVsequence2}) and so $N_p\le L_p$ is equivalent to $\frac{\nu_p}{p}\sum_{k\ge p}\frac{1}{\nu_k}\le 1$ for all such $p$.
\qed\enddemo

Note that \eqref{Lbigequ} can be satisfied for sequences violating \hyperlink{gamma1}{$(\gamma_1)$}. In this case $p\mapsto\frac{\nu_p}{p}\sum_{k\ge p}\frac{1}{\nu_k}$ has to be ''irregular'' (strongly oscillating). For such sequences, having \eqref{Lbigequ} but violating \hyperlink{gamma1}{$(\gamma_1)$}, Lemmas \ref{gamma1} and \ref{Lbig} yield:
\begin{itemize}
\item[$(i)$] $N\hyperlink{approx}{\approx}L$ fails, which means that $N\hyperlink{preceq}{\preceq}L$ is violated but, on the other hand,
\item[$(ii)$] $L$ is ''relatively close'' to $N$ in the sense that even $L_p\ge N_p$ for infinitely many $p$. In particular, in this case $L$ cannot be strictly smaller than $N$, i.e. $\lim_{p\rightarrow+\infty}\left(\frac{L_p}{N_p}\right)^{1/p}=0$ fails.

\item[$(iii)$] Consequently, the failure of $\mathcal{B}$ from being surjective can be viewed as ''relatively small'' within the (mixed) ultradifferentiable weight sequence setting.
\end{itemize}

\subsection{Comparison of the optimal sequences $L$ and $S^N$}\label{descendentcompare}
In this section we will compare the optimal sequence $L$ from the previous section with the descendant $S^N$. The aim is to show that for many cases there is no difference between both optimal sequences which is an advantage because the descendant automatically has good regularity properties and its definition is not as technical as that of $L$.\vspace{6pt}

First, we introduce the following conditions for given nonquasianalytic $N\in\hyperlink{LCset}{\mathcal{LC}}$:
\begin{equation}\label{nongamma2}
\liminf_{p\rightarrow+\infty}\frac{\nu_p}{p}\sum_{k\ge 2p}\frac{1}{\nu_k}>0,
\end{equation}
and
\begin{equation}\label{nongamma2weak}
\liminf_{p\rightarrow+\infty}\frac{\nu_p}{p}\sum_{k\ge p}\frac{1}{\nu_k}>0.
\end{equation}

Note that \eqref{nongamma2} implies \eqref{nongamma2weak} and \eqref{nongamma2} precludes condition $(\gamma_2)$ of \cite{petzsche} with the choice $k=2$ there, while \eqref{nongamma2weak} precludes $(\gamma_2)$ with $k=1$.

Recall that for any given $N\in\hyperlink{LCset}{\mathcal{LC}}$ condition \hyperlink{mg}{$(\on{mg})$} is equivalent to $\sup_{p\in\NN}\frac{\nu_{2p}}{\nu_p}<+\infty$ (see e.g. \cite[Lemma 2.2]{whitneyextensionweightmatrix}). Using this we can prove the following:

\begin{lemma}\label{optimalSV1}
Let $N\in\hyperlink{LCset}{\mathcal{LC}}$ be nonquasianalytic. If $N$ satisfies \hyperlink{mg}{$(\on{mg})$}, then \eqref{nongamma2} holds true, but the converse fails in general.
\end{lemma}

\demo{Proof}
We prove that $\liminf_{p\rightarrow+\infty}\frac{\nu_p}{p}\sum_{k\ge 2p}\frac{1}{\nu_k}=0$ precludes \hyperlink{mg}{$(\on{mg})$}. (If $\liminf$ is replaced by $\lim$ this implication follows by combining \cite[Prop. 1.1. $(b)$, Prop. 1.6 $(a)$]{petzsche} and the comments between Example 1.7 and Example 1.8 in \cite{petzsche}.)\vspace{6pt}

If $\liminf_{p\rightarrow+\infty}\frac{\nu_p}{p}\sum_{k\ge 2p}\frac{1}{\nu_k}=0$, then
$$\forall\;\varepsilon\le 1\;\exists\;p_{\varepsilon}\in\NN_{>0}:\;\;\;\frac{\nu_{p_{\varepsilon}}}{p_{\varepsilon}}\sum_{k\ge 2p_{\varepsilon}}\frac{1}{\nu_k}\le\varepsilon,$$
and so $\sum_{k\ge 2p_{\varepsilon}}\frac{1}{\nu_k}\le\frac{p_{\varepsilon}\varepsilon}{\nu_{p_{\varepsilon}}}$. Since
$$\sum_{k\ge 2p_{\varepsilon}}\frac{1}{\nu_k}=\sum_{k\ge 4p_{\varepsilon}}\frac{1}{\nu_k}+\sum_{k=2p_{\varepsilon}}^{4p_{\varepsilon}-1}\frac{1}{\nu_k}\ge\sum_{k\ge 4p_{\varepsilon}}\frac{1}{\nu_k}+\frac{2p_{\varepsilon}}{\nu_{4p_{\varepsilon}}}\ge\frac{2p_{\varepsilon}}{\nu_{4p_{\varepsilon}}}$$ we arrive at $\frac{\nu_{4p_{\varepsilon}}}{\nu_{p_{\varepsilon}}}\ge\frac{2}{\varepsilon}$. But this contradicts $\sup_{p\in\NN}\frac{\nu_{2p}}{\nu_p}<+\infty$ as $\varepsilon\rightarrow 0$, hence \hyperlink{mg}{$(\on{mg})$} cannot hold true.\vspace{6pt}

The converse fails in general, see \cite[Example 1]{mixedsectorialextensions}.
\qed\enddemo

These arising conditions are related to the technical assumption \hyperlink{mg}{$(\on{mg})$} for the descendant.

\begin{remark}\label{desremark1}
Let $N\in\hyperlink{LCset}{\mathcal{LC}}$ be nonquasianalytic.
\begin{itemize}

\item[$(a)$] In \cite[Lemma 6]{mixedsectorialextensions} a precise characterization has been given when $S^N$ satisfies \hyperlink{mg}{$(\on{mg})$}.

\item[$(b)$] There it has also been shown that \eqref{nongamma2} for $N$ implies \hyperlink{mg}{$(\on{mg})$} for $S^N$.

\item[$(c)$] However, \cite[Example 1]{mixedsectorialextensions} provides an example for a sequence $N$ not satisfying \hyperlink{mg}{$(\on{mg})$} but such that \eqref{nongamma2} holds true (and hence \hyperlink{mg}{$(\on{mg})$} is valid for $S^N$).
\end{itemize}
\end{remark}

Using this observation, in the next result we prove that for many sequences there is no difference between the optimal but technical sequence(s) $L^s$ introduced in this paper and the known descendant $S^N$.

\begin{theorem}\label{optimalSV2}
Let $N\in\hyperlink{LCset}{\mathcal{LC}}$ satisfy \eqref{nongamma2} (which implies nonquasianalyticity). Then the following equivalent assertions are satisfied:
\begin{itemize}
\item[$(i)$] For all $s\in\NN_{>0}$ we have $S^N\hyperlink{approx}{\approx}L^s$.

\item[$(ii)$] $S^N$ is \hyperlink{preceq}{$\preceq$}-maximal among all $M\in\mathcal{N}_{\preceq,R}$ resp. $M\in\mathcal{N}_{\preceq,B}$ (and also among all $M\in\mathcal{N}_{\preceq,\mathcal{LC},R}$ resp. $M\in\mathcal{N}_{\preceq,\mathcal{LC},B}$) satisfying $(M,N)_{SV}$.

\item[$(iii)$] The inclusions $\mathcal{B}(\mathcal{D}_{\{N\}}([-1,1]))\supseteq\Lambda_{\{S^N\}}$ resp. $\mathcal{B}(\mathcal{D}_{(N)}([-1,1]))\supseteq\Lambda_{(S^N)}$ are optimal in the ultradifferentiable setting.
\end{itemize}
\end{theorem}

In particular, this statement holds true for all nonquasianalytic $N\in\hyperlink{LCset}{\mathcal{LC}}$ satisfying \hyperlink{mg}{$(\on{mg})$}.

\demo{Proof}
By Remark \ref{desremark} and Lemma \ref{optimalSV} it remains to show $L^s\hyperlink{preceq}{\preceq}S^N$. By \eqref{nongamma2} the descendant $S^N$ satisfies \hyperlink{mg}{$(\on{mg})$} and so, by Stirling's formula, it suffices to show that
$$\exists\;D\ge 1\;\forall\;p\in\NN_{>0}:\;\;\;(L^s_p/p!)^{1/p}\le D\frac{\sigma_p^N}{p}.$$
Then \eqref{optimalSVsequence} implies for $p\in\NN_{>0}$
$$(L^s_p/p^p)^{1/p}=s\min_{0\le j\le p-1}\left(\frac{1}{\sum_{k\ge p}\frac{1}{\nu_k}}\right)^{1-j/p}p^{-j/p}(N_j)^{1/p}\le s\frac{1}{\sum_{k\ge p}\frac{1}{\nu_k}},$$
by choosing $j=0$. By Stirling's formula it is enough to prove
\begin{align*}
&s\frac{1}{\sum_{k\ge p}\frac{1}{\nu_k}}\le D\frac{\sigma_p^N}{p}=D\frac{\tau_1}{\tau_p}
\\&
\Leftrightarrow\frac{p}{\nu_p}+\sum_{k\ge p}\frac{1}{\nu_k}\le\frac{D\tau_1}{s}\sum_{k\ge p}\frac{1}{\nu_k}\Leftrightarrow\frac{p}{\nu_p}\le\left(\frac{D\tau_1}{s}-1\right)\sum_{k\ge p}\frac{1}{\nu_k}
\\&
\Leftrightarrow\frac{1}{D\tau_1/s-1}\le\frac{\nu_p}{p}\sum_{k\ge p}\frac{1}{\nu_k}.
\end{align*}
The last equivalence holds true if and only if $D\tau_1/s-1>0$ and so if and only if $D>s/\tau_1$. Finally, \eqref{nongamma2} implies \eqref{nongamma2weak}, thus by choosing $D\ge 1$ sufficiently large (depending on given $N$ and $s$) we are done.
\qed\enddemo

{\itshape Conclusion:}

Remark \ref{desremark}, Lemma \ref{gamma1} and Theorem \ref{optimalSV2} show that we have
$$\forall\;s\in\NN_{>0}\;\forall\;C\ge 1:\;\;\;S^N\hyperlink{approx}{\approx}L^{s,C}$$
whenever
$$\liminf_{p\rightarrow+\infty}\frac{\nu_p}{p}\sum_{k\ge 2p}\frac{1}{\nu_k}>0,\hspace{20pt}\text{or}\hspace{20pt}\sup_{p\in\NN_{>0}}\frac{\nu_p}{p}\sum_{k\ge p}\frac{1}{\nu_k}<+\infty,$$
and in the latter case we even have $S^N\hyperlink{approx}{\approx}L^{s,C}\hyperlink{approx}{\approx}N$.

\subsection{The ramified case}\label{ramisection}
For completeness, we now give some comments on the ramified framework as well, with a ramification parameter $r>0$. This setting and the corresponding conditions become meaningful when treating (nonstandard) ultradifferentiable ramification classes from \cite{Schmetsvaldivia00} which are needed for the study of the surjectivity of the (asymptotic) Borel map and to prove (mixed) extension results in the {\itshape ultraholomorphic setting.} The parameter $r$ is then used to introduce the (mixed) growth index $\gamma(\cdot)$ which measures the opening of the sector under consideration. For more details we refer to \cite{mixedramisurj}, \cite{injsurj} and \cite{mixedsectorialextensions} and the references therein.\vspace{6pt}

We have the crucial mixed conditions
\begin{equation*}\label{SVrami}
\exists\;s\in\NN_{>0}:\;\;\sup_{p\in\NN_{>0}}\frac{(\lambda_{p,s}^{M,N})^{1/r}}{p}\sum_{k\ge p}\left(\frac{1}{\nu_k}\right)^{1/r}<+\infty,
\end{equation*}
denoted by \hypertarget{SV}{$(M,N)_{SV_r}$}, and
\begin{equation*}\label{gamma1rami}
\sup_{p\in\NN_{>0}}\frac{(\mu_p)^{1/r}}{p}\sum_{k\ge p}\left(\frac{1}{\nu_k}\right)^{1/r}<+\infty,
\end{equation*}
denoted by \hypertarget{gammarmix}{$(M,N)_{\gamma_r}$}. For each $h>0$ and $r\in\NN_{>0}$ one considers the Banach space
\begin{align*}
\mathcal{D}_{r,M,h}([-1,1]):=&\Big\{f\in\mathcal{E}(\RR,\CC):\;\supp(f)\subseteq[-1,1],
\\&
f^{(rp+j)}(0)=0\;\forall\;p\in\NN,\;\forall\;j\in\{1,\dots,r-1\},\;\;\sup_{p\in\NN, x\in\RR}\frac{|f^{(rp)}(x)|}{h^pM_p}<+\infty\Big\},
\end{align*}
and the $r$-ramified ultradifferentiable test function class of Roumieu-type by
\begin{equation*}\label{Roumieu-LB-space}
\mathcal{D}_{r,\{M\}}([-1,1]):=\underset{h>0}{\varinjlim}\;\mathcal{D}_{r,M,h}([-1,1]),
\end{equation*}
respectively the Beurling-type by
\begin{equation*}\label{Beurling-Frechet-space}
\mathcal{D}_{r,(M)}([-1,1]):=\underset{h>0}{\varprojlim}\;\mathcal{D}_{r,M,h}([-1,1]).
\end{equation*}
Finally, the {\itshape ramified Borel map} (at $0$) is given by
$$\mathcal{B}^r:\mathcal{D}_{r,[M]}([-1,1])\longrightarrow\Lambda_{[M]},\hspace{20pt}\mathcal{B}^r(f):=(f^{(rp)}(0))_{p\in\NN}.$$

When applying the above definitions and techniques to $N^{1/r}$ (instead of $N$), all results can be transferred in an obvious and straightforward way to this framework. The optimal sequence for $(M,N)_{SV_r}$ is given by $L^{r;s}:=\left(L^{N^{1/r};s}\right)^r$, with $L^{N;s}:=L^s$, whereas the optimal sequence for $(M,N)_{\gamma_r}$ is $M^{N,r}:=(S^{N,r})^r$, with $S^{N,r}$ denoting the descendant of $N^{1/r}$.

\section{A (Counter)-example on the optimal sequence}\label{counterex}
The aim is to show that in general $S^N\hyperlink{approx}{\approx}L^{s,C}$ will fail, by constructing an appropriate $N$. Thus $L$ is in general more optimal than $S^N$. The previous conclusion and Lemma \ref{Lbig} indicate that the oscillation or irregularity of $p\mapsto\frac{\nu_p}{p}\sum_{k\ge p}\frac{1}{\nu_k}$ has to show up in order to destroy the equivalence. First we record the following easy observation:

\begin{lemma}
Let $N\in\hyperlink{LCset}{\mathcal{LC}}$. If $S^N\hyperlink{approx}{\approx}L^{s,C}$ is violated, then \hyperlink{gamma1}{$(\gamma_1)$} for $N$ does not hold and $L^{s,C}\hyperlink{approx}{\approx}N$ fails.
\end{lemma}

\demo{Proof}
We can assume that $N$ is nonquasianalytic, otherwise the definition of $L^{s,C}$ would not make sense (see also Remark \ref{lquasianalytic}). If now $N$ satisfies \hyperlink{gamma1}{$(\gamma_1)$}, then $S^N\hyperlink{approx}{\approx}N$ follows and so $S^N\hyperlink{preceq}{\preceq}L^{s,C}\hyperlink{preceq}{\preceq}N$ implies $S^N\hyperlink{approx}{\approx}L^{s,C}$.

Hence, if $S^N\hyperlink{approx}{\approx}L^{s,C}$ is violated then \hyperlink{gamma1}{$(\gamma_1)$} has to fail and, by Lemma \ref{gamma1}, $L^{s,C}\hyperlink{approx}{\approx}N$ cannot be valid.
\qed\enddemo

To come up with an example, we need some preparation and have to recall some statements.\vspace{6pt}

The starting idea has been suggested by Armin Rainer. Since $s^N:=(S^N_p/p!)_{p\in\NN}\in\hyperlink{LCset}{\mathcal{LC}}$, i.e. the descendant $S^N$ is normalized and satisfies \hyperlink{slc}{$(\text{slc})$}, the ultradifferentiable classes $\mathcal{E}_{\{S^N\}}$ and $\mathcal{E}_{(S^N)}$ are closed under composition and have important stability properties, see \cite{compositionpaper} and \cite{characterizationstabilitypaper}.

Since $s^N$ is normalized and log-convex, the mapping $p\mapsto(s^N_p)^{1/p}$ is nondecreasing and consequently for each sequence $M\in\RR_{>0}^{\NN}$ with $S^N\hyperlink{approx}{\approx}M$, which is equivalent to $s^N\hyperlink{approx}{\approx}m$, we get $$\exists\;C\ge 1\;\forall\;1\le j\le k:\;\;\;\frac{1}{C}(m_j)^{1/j}\le(s^N_j)^{1/j}\le(s^N_k)^{1/k}\le C(m_k)^{1/k}.$$
Thus each $M$ which is equivalent to $S^N$ has to satisfy
\begin{equation}\label{almostincreasing}
\exists\;C\ge 1\;\forall\;1\le p\le q:\;\;\;(m_p)^{1/p}\le C(m_q)^{1/q},
\end{equation}
i.e. the sequence $(m_p^{1/p})_{p\ge 1}$ has to be {\itshape almost increasing.}

Applying this information to $M\equiv L^s$, the equivalence $S^N\hyperlink{approx}{\approx}L^s$ implies (with Stirling's formula) that the sequence $((L^s_p/p^p)^{1/p})_{p\ge 1}$ has to be almost increasing. Since $L^s\hyperlink{approx}{\approx}L^t$ for all $s,t\in\NN_{>0}$ we see that in this situation $((L_p/p^p)^{1/p})_{p\ge 1}$ has to be almost increasing and
$$(L_p/p^p)^{1/p}=\min_{0\le j\le p-1}\left(\frac{1}{\sum_{k\ge p}\frac{1}{\nu_k}}\right)^{1-j/p}p^{-j/p}(N_j)^{1/p}.$$

Given $p\in\NN_{>0}$ and $0\le j<j+1\le p-1$ we get
\begin{align*}
&\left(\frac{1}{\sum_{k\ge p}\frac{1}{\nu_k}}\right)^{1-j/p}p^{-j/p}(N_j)^{1/p}\ge\left(\frac{1}{\sum_{k\ge p}\frac{1}{\nu_k}}\right)^{1-(j+1)/p}p^{-(j+1)/p}(N_{j+1})^{1/p}
\\&
\Leftrightarrow 1\ge\left(\frac{\nu_{j+1}}{p}\right)^{1/p}\left(\sum_{k\ge p}\frac{1}{\nu_k}\right)^{1/p}\Leftrightarrow 1\ge\frac{\nu_{j+1}}{p}\sum_{k\ge p}\frac{1}{\nu_k}.
\end{align*}
Thus, if
\begin{equation}\label{minlastpossnec}
1\ge\frac{\nu_{p-1}}{p}\sum_{k\ge p}\frac{1}{\nu_k},
\end{equation}
then by log-convexity (i.e. $j\mapsto\nu_j$ is nondecreasing) we get
\begin{equation}\label{minlastposs}
(L_p/p^p)^{1/p}=\min_{0\le j\le p-1}\left(\frac{1}{\sum_{k\ge p}\frac{1}{\nu_k}}\right)^{1-j/p}p^{-j/p}(N_j)^{1/p}\underbrace{=}_{j=p-1}\left(\frac{1}{\sum_{k\ge p}\frac{1}{\nu_k}}\right)^{1/p}\frac{p^{1/p}}{p}(N_{p-1})^{1/p}.
\end{equation}
Hence, in order to show that \eqref{almostincreasing} is violated for $L$, it suffices to prove that for any $C\ge 1$ large we can find an integer $p$ satisfying \eqref{minlastpossnec} and some other integer $q>p$ such that
\begin{equation}\label{violater}
(L_p/p^p)^{1/p}=\left(\frac{1}{\sum_{k\ge p}\frac{1}{\nu_k}}\right)^{1/p}\frac{p^{1/p}}{p}(N_{p-1})^{1/p}\ge C\left(\frac{1}{\sum_{k\ge q}\frac{1}{\nu_k}}\right)^{1/q}\frac{q^{1/q}}{q}(N_{q-1})^{1/q}(\ge C(L_q/q^q)^{1/q}).
\end{equation}
If
\begin{equation}\label{liminfsmaller1}
\liminf_{p\rightarrow+\infty}\frac{\nu_{p}}{p}\sum_{k\ge p}\frac{1}{\nu_k}<1,
\end{equation}
then \eqref{minlastpossnec} and so \eqref{minlastposs} holds true for infinitely many $p$ and for each such $p$ we have $p^{1/p}\frac{1}{\left(\sum_{k\ge p}\frac{1}{\nu_k}\right)^{1/p}}\ge(\nu_{p})^{1/p}$. Hence, the left-hand side in \eqref{minlastposs} yields, for such $p(\ge 2)$
$$(L_p/p^p)^{1/p}=\left(\frac{1}{\sum_{k\ge p}\frac{1}{\nu_k}}\right)^{1/p}\frac{p^{1/p}}{p}(N_{p-1})^{1/p}\ge\frac{(\nu_{p})^{1/p}}{p}(N_{p-1})^{1/p}=\frac{(N_{p})^{1/p}}{p}.$$

\eqref{liminfsmaller1} is connected to condition \eqref{nongamma2} as follows:

\begin{lemma}
Let $N\in\hyperlink{LCset}{\mathcal{LC}}$ be nonquasianalytic. Then
$$\liminf_{p\rightarrow+\infty}\frac{\nu_p}{p}\sum_{k\ge 2p}\frac{1}{\nu_k}=0,$$
i.e. $\neg$\eqref{nongamma2}, implies $\liminf_{p\rightarrow+\infty}\frac{\nu_p}{p}\sum_{k\ge p}\frac{1}{\nu_k}\le 1$.
\end{lemma}

\demo{Proof}
If $\liminf_{p\rightarrow+\infty}\frac{\nu_p}{p}\sum_{k\ge p}\frac{1}{\nu_k}>1$, then
$$\exists\;\varepsilon>0\;\exists\;p_{\varepsilon}\in\NN\;\forall\;p\ge p_{\varepsilon}:\;\;\;\frac{\nu_p}{p}\sum_{k\ge p}\frac{1}{\nu_k}\ge 1+\varepsilon.$$
Hence, for all $p\ge p_{\varepsilon}$
$$\frac{\nu_p}{p}\sum_{k\ge p}\frac{1}{\nu_k}=\frac{\nu_p}{p}\sum_{k\ge 2p}\frac{1}{\nu_k}+\frac{\nu_p}{p}\sum_{k=p}^{2p-1}\frac{1}{\nu_k}\ge 1+\varepsilon\Longleftrightarrow\frac{\nu_p}{p}\sum_{k\ge 2p}\frac{1}{\nu_k}\ge 1+\varepsilon-\frac{\nu_p}{p}\sum_{k=p}^{2p-1}\frac{1}{\nu_k},$$
and
$$1+\varepsilon-\frac{\nu_p}{p}\sum_{k=p}^{2p-1}\frac{1}{\nu_k}\ge\varepsilon\Longleftrightarrow 1\ge\frac{\nu_p}{p}\sum_{k=p}^{2p-1}\frac{1}{\nu_k},$$
which holds true because $\frac{\nu_p}{p}\sum_{k=p}^{2p-1}\frac{1}{\nu_k}\le\frac{\nu_p}{p}\frac{p}{\nu_p}=1$. Consequently, we get $\frac{\nu_p}{p}\sum_{k\ge 2p}\frac{1}{\nu_k}\ge\varepsilon>0$ for all $p\ge p_{\varepsilon}$ which proves \eqref{nongamma2}.
\qed\enddemo

Similarly, when
$$\limsup_{q\rightarrow+\infty}\frac{\nu_q}{q}\sum_{k\ge q}\frac{1}{\nu_k}>A\ge 1,$$
then for infinitely many integers $q$,
$$q^{1/q}\frac{1}{\left(\sum_{k\ge q}\frac{1}{\nu_k}\right)^{1/q}}<\frac{(\nu_q)^{1/q}}{A^{1/q}},$$
and for each such $q$ we get
$$(L_q/q^q)^{1/q}\le\left(\frac{1}{\sum_{k\ge q}\frac{1}{\nu_k}}\right)^{1/q}\frac{q^{1/q}}{q}(N_{q-1})^{1/q}\le\frac{1}{A^{1/q}}\frac{(\nu_q)^{1/q}}{q}(N_{q-1})^{1/q}=\frac{1}{A^{1/q}}\frac{(N_q)^{1/q}}{q}.$$
Summarizing, by recalling that by Stirling's formula $p!^{1/p}$ and $p$ grow similarly up to a constant, we have shown the following:

\begin{lemma}\label{maximalviolatedlemma}
Let $N\in\hyperlink{LCset}{\mathcal{LC}}$ be nonquasianalytic. Assume that for all $i\in\NN_{>0}$ (or more generally for all $i\ge i_0\in\NN_{>0}$) there are $p_i,q_i\in\NN_{>0}$ with $p_{i+1}>q_i>p_i$, and constants $A_i\ge 1$, $C_i>0$, with $\lim_{i\rightarrow+\infty}\frac{(A_i)^{1/q_i}}{C_i}=+\infty$, such that
\begin{itemize}
\item[$(I)$] $\frac{\nu_{p_i}}{p_i}\sum_{k\ge p_i}\frac{1}{\nu_k}<1$,

\item[$(II)$] $\frac{\nu_{q_i}}{q_i}\sum_{k\ge q_i}\frac{1}{\nu_k}>A_i$,

\item[$(III)$] $C_i(n_{p_i})^{1/p_i}\ge(n_{q_i})^{1/q_i}$.
\end{itemize}

Then $((L_p/p^p)^{1/p})_{p\ge 1}$ is not almost increasing, hence $S^N\hyperlink{approx}{\approx}L^s$ cannot hold true.

\end{lemma}

We now construct $N$ satisfying the requirements of Lemma \ref{maximalviolatedlemma}. $N$ will be defined in terms of the sequence of quotients $(\nu_j)_{j\ge 1}$ and for this let $(a_j)_{j\in\NN_{>0}}$, $(b_j)_{j\in\NN_{>0}}$ be strictly increasing sequences such that $a_j>1$ for all $j\ge 1$ and $\lim_{j\rightarrow+\infty}a_j=\lim_{j\rightarrow+\infty}b_j=+\infty$.

The required sequences $(p_i)_i$ and $(q_i)_i$ (in $\NN_{>0}$) will be defined iteratively with $p_{i+1}>q_i>p_i$ and increasing fast enough to guarantee at least (for given $p_i$)
\begin{equation}\label{qipirequirement}
p_1>2,\;\;\;\forall\;i\in\NN_{>0}:\;\;\;q_i\ge 2p_i-1\Leftrightarrow\frac{q_i-p_i}{p_i-1}\ge 1.
\end{equation}
For the sequence $(a_j)_j$ we require that
\begin{equation}\label{airequirement}
\forall\;i\in\NN_{>0}:\;\;\;\frac{p_{i+1}-2}{p_i-2}<\frac{q_i-p_i}{p_i-2}a_i\Longrightarrow\frac{p_{i+1}-q_i}{p_i-2}<\frac{q_i-p_i}{p_i-2}a_i,
\end{equation}
i.e. $a_i$ depends on $p_i$, $q_i$ and $p_{i+1}$.

We make the ansatz
\begin{equation}\label{ansatzpiqi}
\frac{N_{q_i}}{q_i!}=n_{q_i}=(n_{p_i})^{q_i/p_i}C_i^{q_i}=\left(\frac{N_{p_i}}{p_i!}\right)^{q_i/p_i}C_i^{q_i}\Longleftrightarrow N_{q_i}=q_i!\left(\frac{N_{p_i}}{p_i!}\right)^{q_i/p_i}C_i^{q_i},
\end{equation}
with $C_i\ge 1$, $C_i\rightarrow+\infty$ as $i\rightarrow+\infty$ such that
\begin{equation}\label{sequencecrequirebasic}
C_1\ge\frac{(p_1!)^{(q_1-p_1)/(p_1q_1)}}{((p_1+1)\cdots q_1)^{1/q_1}},
\end{equation}
and such that (we have $p_i\ge p_2>2$, $i\ge 2$)
\begin{equation}\label{sequencecrequire}
C_{i+1}\ge C_i^{\frac{q_i(q_{i+1}-p_{i+1})}{q_{i+1}(q_i-p_i)}}\left(a_i\frac{(p_{i+1}!)^{1/p_{i+1}}(N_{p_i})^{1/p_i}}{(p_i!)^{1/p_i}(N_{p_{i+1}})^{1/p_{i+1}}}2^{i+2}\frac{q_{i+1}-p_{i+1}}{p_{i+1}-2}\right)^{(q_{i+1}-p_{i+1})/q_{i+1}},\;\;\;i\ge 1.
\end{equation}
Consequently, by taking into account \eqref{qipirequirement}, we get
\begin{equation}\label{sequencecrequire1}
C_{i+1}\ge C_i^{\frac{q_i(q_{i+1}-p_{i+1})}{q_{i+1}(q_i-p_i)}}\left(a_i\frac{(p_{i+1}!)^{1/p_{i+1}}(N_{p_i})^{1/p_i}}{(p_i!)^{1/p_i}(N_{p_{i+1}})^{1/p_{i+1}}}\right)^{(q_{i+1}-p_{i+1})/q_{i+1}},\;\;\;i\ge 1.
\end{equation}

Moreover we set
\begin{equation}\label{sequenceachoice}
A_i:=(C_i\cdot b_i)^{q_i},\;\;\;i\in\NN_{>0},
\end{equation}
which yields $\lim_{i\rightarrow+\infty}\frac{(A_i)^{1/q_i}}{C_i}=+\infty$.

Then, by \eqref{ansatzpiqi}, $$\nu_{p_i+1}\cdots\nu_{q_i}=\frac{N_{q_i}}{N_{p_i}}=C_i^{q_i}\frac{q_i!(N_{p_i})^{q_i/p_i}}{(p_i!)^{q_i/p_i}N_{p_i}}=C_i^{q_i}\frac{q_i!}{p_i!}(p_i!)^{(p_i-q_i)/p_i}(N_{p_i})^{(q_i-p_i)/p_i},$$ and so, when taking $\nu_{p_i+1}=\dots=\nu_{q_i}$, we have to put
\begin{equation}\label{nuchoic1}
\nu_{p_i+1}=\dots=\nu_{q_i}:=C_i^{q_i/(q_i-p_i)}(N_{p_i})^{1/p_i}\frac{1}{(p_i!)^{1/p_i}}((p_i+1)\cdots q_i)^{1/(q_i-p_i)},\;\;\;i\in\NN_{>0}.
\end{equation}
Moreover, we set
\begin{equation}\label{nuchoic2}
\nu_{q_i+1}=\dots=\nu_{p_{i+1}}:=a_i\cdot\nu_{q_i}=a_i\cdot\nu_{p_i+1},\;\;\;i\in\NN_{>0},
\end{equation}
and finally, in order to complete the definition of $N$:
\begin{equation}\label{nuchoic3}
\nu_0=\dots=\nu_{p_1}:=1.
\end{equation}
Note that by \eqref{sequencecrequire} the choice of $C_{i+1}$ only depends on given values $C_i$, $q_i$, $q_{i+1}$, $p_{i}$, $p_{i+1}$, on $a_i$ (related to $p_i$, $q_i$ and $p_{i+1}$ via \eqref{airequirement}) and finally on $N_{p_i}$ and $N_{p_{i+1}}$, involving again only terms depending on $a_i$, $C_i$, $q_i$ and $p_i$. Hence this choice of $C_{i+1}$ is then possible and used to determine $N_{q_{i+1}}$ via \eqref{ansatzpiqi}, hence $\nu_{p_{i+1}+1},\dots,\nu_{q_{i+1}},\nu_{q_{i+1}+1},\dots,\nu_{q_{i+2}}$ via the above definitions.\vspace{6pt}

{\itshape Claim:} $N\in\hyperlink{LCset}{\mathcal{LC}}$ holds true. Normalization, i.e. $1=N_0\le N_1$, follows because $1=\nu_0=\nu_1$ (recall $p_1>1$). For $k\mapsto\nu_k$ to be nondecreasing, first by \eqref{nuchoic1} and \eqref{nuchoic3}, we need to check $\nu_{p_1}\le\nu_{p_1+1}$, so
\begin{align*}
1&\le C_1^{q_1/(q_1-p_1)}(N_{p_1})^{1/p_1}\frac{1}{(p_1!)^{1/p_1}}((p_1+1)\cdots q_1)^{1/(q_1-p_1)}
\\&
=C_1^{q_1/(q_1-p_1)}\frac{1}{(p_1!)^{1/p_1}}((p_1+1)\cdots q_1)^{1/(q_1-p_1)},
\end{align*}
because $N_{p_1}=1$ by \eqref{nuchoic3}. This is valid by the choice of $C_1$ in \eqref{sequencecrequirebasic}.

Second, by \eqref{nuchoic1} and \eqref{nuchoic2} and because $a_i>1$ it suffices to check $\nu_{p_{i+1}+1}\ge\nu_{p_{i+1}}$, $i\ge 1$, which is equivalent to
\begin{align*}
&C_{i+1}^{q_{i+1}/(q_{i+1}-p_{i+1})}(N_{p_{i+1}})^{1/p_{i+1}}\frac{1}{(p_{i+1}!)^{1/p_{i+1}}}((p_{i+1}+1)\cdots q_{i+1})^{1/(q_{i+1}-p_{i+1})}
\\&
\ge\nu_{q_i+1}=a_i\nu_{p_i+1}=a_iC_i^{q_i/(q_i-p_i)}(N_{p_i})^{1/p_i}\frac{1}{(p_i!)^{1/p_i}}((p_i+1)\cdots q_i)^{1/(q_i-p_i)}.
\end{align*}
We get
$$((p_i+1)\cdots q_i)^{1/(q_i-p_i)}\le q_i^{(q_i-p_i)/(q_i-p_i)}=q_i\le p_{i+1}+1\le((p_{i+1}+1)\cdots q_{i+1})^{1/(q_{i+1}-p_{i+1})}$$ and by taking into account \eqref{sequencecrequire1} we have shown this inequality. Again, by \eqref{nuchoic1} and \eqref{nuchoic2} and since $\lim_{j\rightarrow+\infty}a_j=+\infty$ we obtain $\lim_{p\rightarrow+\infty}\nu_p=+\infty$, hence $\lim_{p\rightarrow+\infty}(N_p)^{1/p}=+\infty$ as well.\vspace{6pt}

{\itshape Claim:} Requirement $(III)$ holds true (with equality for all $i\in\NN_{>0}$). This is immediate by \eqref{ansatzpiqi}.\vspace{6pt}

{\itshape Claim:} Requirement $(II)$ holds true. First, we get
$$\frac{\nu_{q_i}}{q_i}\sum_{k\ge q_i}\frac{1}{\nu_k}=\frac{1}{q_i}+\frac{\nu_{q_i}}{q_i}\sum_{k\ge q_i+1}\frac{1}{\nu_k}>A_i\Leftrightarrow\sum_{k\ge q_i+1}\frac{1}{\nu_k}>\frac{A_iq_i-1}{\nu_{q_i}}.$$
By \eqref{nuchoic2} we see that $\sum_{k\ge q_i+1}\frac{1}{\nu_k}\ge\frac{p_{i+1}-q_i}{\nu_{q_i+1}}$ holds and because all the arising further summands are positive it suffices to show (see \eqref{sequenceachoice}) that $$\forall\;i\in\NN_{>0}:\;\;\;\frac{p_{i+1}-q_i}{\nu_{q_i+1}}>\frac{A_iq_i-1}{\nu_{q_i}}=\frac{(b_iC_i)^{q_i}q_i-1}{\nu_{q_i}},$$
which is equivalent to
$$\forall\;i\in\NN_{>0}:\;\;\;p_{i+1}>\frac{\nu_{q_i+1}((b_iC_i)^{q_i}q_i-1)}{\nu_{q_i}}+1.$$
This can be done by choosing $p_{i+1}$ large enough (depending on given $q_i$, $b_i$ and $C_i$). For this recall that $C_i$ only depends on given $C_{i-1}$, $q_{i-1}$, $q_i$, $p_{i-1}$, $p_i$, on $a_{i-1}$ (related to $p_i$, $q_i$ and $p_i$ via \eqref{airequirement}) and on $N_{p_{i-1}}$ and $N_{p_i}$ via \eqref{sequencecrequire}. Finally, $\nu_{q_i+1}$ also only depends on these values, see \eqref{nuchoic1} and \eqref{nuchoic2}.\vspace{6pt}

{\itshape Claim:} Requirement $(I)$ holds true. We see that
$$\frac{\nu_{p_i}}{p_i}\sum_{k\ge p_i}\frac{1}{\nu_k}=\frac{1}{p_i}+\frac{\nu_{p_i}}{p_i}\sum_{k\ge p_i+1}\frac{1}{\nu_k}\le\frac{1}{2}$$
and this is equivalent to $\sum_{k\ge p_i+1}\frac{1}{\nu_k}\le\frac{p_i-2}{2\nu_{p_i}}$.

For the series on the left-hand side we see by \eqref{nuchoic1} and \eqref{nuchoic2} that
$$\forall\;k\ge i\ge 1:\;\;\;\sum_{k\ge p_i+1}\frac{1}{\nu_k}=\sum_{k\ge i}\frac{q_k-p_k}{\nu_{p_k+1}}+\sum_{k\ge i}\frac{p_{k+1}-q_k}{\nu_{q_k+1}}=:\sum_{k\ge i}\alpha_k+\sum_{k\ge i}\beta_k.$$
Thus, in order to conclude it suffices to verify
\begin{equation}\label{alphabetaverify}
\forall\;k\ge i\ge 2:\;\;\;\alpha_k<\frac{1}{2^{k+1}}\frac{p_i-2}{\nu_{p_i}},\hspace{30pt}\beta_k<\frac{1}{2^{k+1}}\frac{p_i-2}{\nu_{p_i}},
\end{equation}
because then
$$\forall\;i\ge 2:\;\;\;\sum_{k\ge p_i+1}\frac{1}{\nu_k}=\sum_{k\ge i}\alpha_k+\sum_{k\ge i}\beta_k<2\frac{p_i-2}{\nu_{p_i}}\sum_{k\ge i}\frac{1}{2^{k+1}}\le\frac{1}{2}\frac{p_i-2}{\nu_{p_i}}.$$
\vspace{6pt}

First let $k=i\ge 2$. Then (recall that $p_i\ge p_2>2$)
\begin{align*}
&\alpha_i<\frac{1}{2^{i+1}}\frac{p_i-1}{\nu_{p_i}}\Leftrightarrow\frac{q_i-p_i}{\nu_{p_i+1}}<\frac{1}{2^{i+1}}\frac{p_i-2}{\nu_{p_i}}\Leftrightarrow 2^{i+1}\frac{q_i-p_i}{p_i-2}\nu_{p_i}=2^{i+1}\frac{q_i-p_i}{p_i-2}a_{i-1}\nu_{p_{i-1}+1}<\nu_{p_i+1}
\\&
\Leftrightarrow 2^{i+1}\frac{q_i-p_i}{p_i-2}a_{i-1}(C_{i-1})^{q_{i-1}/(q_{i-1}-p_{i-1})}(N_{p_{i-1}})^{1/p_{i-1}}\frac{1}{(p_{i-1}!)^{1/p_{i-1}}}((p_{i-1}+1)\cdots q_{i-1})^{1/(q_{i-1}-p_{i-1})}
\\&
<(C_i)^{q_i/(q_i-p_i)}(N_{p_i})^{1/p_i}\frac{1}{(p_i!)^{1/p_i}}((p_i+1)\cdots q_i)^{1/(q_i-p_i)}.
\end{align*}
This holds true by the choice of $C_i$ in \eqref{sequencecrequire} and because
$$(p_{i-1}+1)\cdots q_{i-1})^{1/(q_{i-1}-p_{i-1})}\le q_{i-1}<p_i+1\le((p_i+1)\cdots q_i)^{1/(q_i-p_i)}.$$
\vspace{6pt}

For $\beta_i$ we have
\begin{align*}
&\beta_i<\frac{1}{2^{i+1}}\frac{p_i-2}{\nu_{p_i}}\Leftrightarrow\frac{p_{i+1}-q_i}{\nu_{q_i+1}}<\frac{1}{2^{i+1}}\frac{p_i-2}{\nu_{p_i}}
\\&
\Leftrightarrow 2^{i+1}\frac{p_{i+1}-q_i}{p_i-2}\nu_{p_i}=2^{i+1}\frac{p_{i+1}-q_i}{p_i-2}a_{i-1}\nu_{p_{i-1}+1}<\nu_{q_i+1}=a_i\nu_{p_i+1}
\\&
\Leftrightarrow 2^{i+1}\frac{p_{i+1}-q_i}{p_i-2}a_{i-1}C_{i-1}^{q_{i-1}/(q_{i-1}-p_{i-1})}(N_{p_{i-1}})^{1/p_{i-1}}\frac{1}{(p_{i-1}!)^{1/p_{i-1}}}((p_{i-1}+1)\cdots q_{i-1})^{1/(q_{i-1}-p_{i-1})}
\\&
<a_iC_i^{q_i/(q_i-p_i)}(N_{p_i})^{1/p_i}\frac{1}{(p_i!)^{1/p_i}}((p_i+1)\cdots q_i)^{1/(q_i-p_i)},
\end{align*}
which holds because
$$((p_{i-1}+1)\cdots q_{i-1})^{1/(q_{i-1}-p_{i-1})}\le q_{i-1}<p_i+1\le((p_i+1)\cdots q_i)^{1/(q_i-p_i)},$$
in view of $\frac{p_{i+1}-q_i}{p_i-2}a_{i-1}<\frac{q_i-p_i}{p_i-2}a_{i-1}a_i$ (see \eqref{airequirement}) and finally again by \eqref{sequencecrequire}.\vspace{6pt}

Hence we have checked \eqref{alphabetaverify} for all $k=i\ge 2$ and so, in order to verify \eqref{alphabetaverify}, we now prove, for all $i\ge 2$,
\begin{equation}\label{alphabetaverify1}
\frac{p_{i+1}-2}{\nu_{p_{i+1}}}<\frac{p_i-2}{\nu_{p_i}}\Leftrightarrow\frac{p_{i+1}-2}{p_i-2}\nu_{p_i}<\nu_{p_{i+1}}\Leftrightarrow\frac{p_{i+1}-2}{p_i-2}a_{i-1}\nu_{p_{i-1}+1}<a_i\nu_{p_i+1},
\end{equation}
i.e.
\begin{align*}
&\frac{p_{i+1}-2}{p_i-2}a_{i-1}C_{i-1}^{q_{i-1}/(q_{i-1}-p_{i-1})}(N_{p_{i-1}})^{1/p_{i-1}}\frac{1}{(p_{i-1}!)^{1/p_{i-1}}}((p_{i-1}+1)\cdots q_{i-1})^{1/(q_{i-1}-p_{i-1})}
\\&
<a_iC_i^{q_i/(q_i-p_i)}(N_{p_i})^{1/p_i}\frac{1}{(p_i!)^{1/p_i}}((p_i+1)\cdots q_i)^{1/(q_i-p_i)},
\end{align*}
which follows from $\frac{p_{i+1}-2}{p_i-2}a_{i-1}<\frac{q_i-p_i}{p_i-2}a_{i-1}a_i$ (see \eqref{airequirement}) and once again from \eqref{sequencecrequire}.

\section{Failure of injectivity and surjectivity in the nonquasianalytic setting}\label{failureinjandsurj}
\subsection{Failure of injectivity}\label{failureinj}
The aim of this section is to study the failure for $\mathcal{B}$ to be injective in the nonquasianalytic ultradifferentiable weight sequence setting. We show that this failure is large when measured by vector space dimension.\vspace{6pt}

We call a given smooth function $f$ {\itshape flat} (at $x=0$) if $f^{(k)}(0)=0$ for all $k\in\NN$.\vspace{6pt}

First we recall \cite[Theorem 2.2]{petzsche}, which follows from \cite[Theorem 1.3.5]{hoermander}, for a proof see also \cite[Lemma 5.1.6]{diploma}.

\begin{lemma}\label{Petzsche2}
Let $N\in\hyperlink{LCset}{\mathcal{LC}}$ and assume that $a:=\sum_{j=1}^{+\infty}\frac{1}{\nu_j}<+\infty$, i.e. $N$ is nonquasianalytic. Then there exists a smooth function $\varphi$ with support in $[-a,a]$, such that $0\le\varphi(x)\le 1$ for all $x\in[-a,a]$, and $\varphi^{(j)}(0)=\delta_{j,0}$ (Kronecker delta). Furthermore, $\left\|\varphi^{(j)}\right\|_{\infty}\le 2^jN_j$ for all $j\in\NN$.\vspace{6pt}

So $\varphi$ is a nontrivial function ($\varphi(0)=1$) with compact support and $\varphi\in\mathcal{E}_{\{M\}}(\RR,\CC)$ (take $h=2$).
\end{lemma}

In fact, by inspecting the proof of \cite[Lemma 5.1.6]{diploma} we see that $\supp(\varphi)=[-a,a]$.

We put $\varphi_c(x):=\varphi(cx)$ and $b:=\frac{a}{c}$. So the support of $\varphi_c$ equals $[-\frac{a}{c},\frac{a}{c}]=[-b,b]$ and $\varphi_c^{(j)}(x)=c^j\varphi^{(j)}(cx)$ for all $x\in[-b,b]$, which implies $\varphi_c^{(j)}(0)=0$ for $j\ge 1$ and $\varphi_c(0)=\varphi(0)=1$, i.e. $\varphi^{(j)}_c(0)=\delta_{0,j}$. If $c\ge a$, then clearly $\varphi_c\in\mathcal{D}_{\{N\}}([-1,1])$ (take $h=2c$). On the other hand, given a nonquasianalytic $N$ with $a:=\sum_{j=1}^{+\infty}\frac{1}{\nu_j}<+\infty$, we can replace $N$ by the equivalent sequence $\widetilde{N}=(N_j/a^j)_j$ which satisfies $\sum_{j=1}^{+\infty}\frac{1}{\widetilde{\nu}_j}=1$.\vspace{6pt}

For the Beurling type we recall that in the proof of \cite[Theorem 2.1 $(a)(i)$]{petzsche} even a sequence $(\chi_p)_{p\in\NN}$ of functions with compact support in $\mathcal{E}_{(N)}(\RR,\CC)$ has been constructed which satisfies $\chi^{(j)}_p(0)=\delta_{j,p}$. So all the above holds for the Beurling type classes $\mathcal{D}_{(N)}([-1,1])$ as well by choosing $\varphi:=\chi_0$ and rescaling. Note that in the proof of \cite[Theorem 2.1 $(a)(i)$]{petzsche} Lemma \ref{Petzsche2} has been used for the construction.

\begin{proposition}\label{noninj}
Let $N\in\hyperlink{LCset}{\mathcal{LC}}$ and $M\in\RR_{>0}^{\NN}$ with $M\in\mathcal{N}_{\preceq,R}$ resp. $M\in\mathcal{N}_{\preceq,B}$ and such that $(M,N)_{SV}$ holds true.

Then for any $\mathbf{a}\in\Lambda_{[M]}$ there exist infinitely many functions in $\mathcal{D}_{[N]}([-1,1])$ which are mapped onto $\mathbf{a}$ by $\mathcal{B}$. In fact, we can assume that these functions form an infinite-dimensional affine vector space of dimension $\mathfrak{c}$.
\end{proposition}

\demo{Proof}
By Theorem \ref{SVtheorem}, first we see that $N$ is nonquasianalytic and so $a:=\sum_{j=1}^{+\infty}\frac{1}{\nu_j}<+\infty$, and second $\mathcal{B}(\mathcal{D}_{\{N\}}([-1,1]))\supseteq\Lambda_{\{M\}}$ resp. $\mathcal{B}(\mathcal{D}_{(N)}([-1,1]))\supseteq\Lambda_{(M)}$.

Thus for each given $\mathbf{a}\in\Lambda_{[M]}$ we can find some $f_{\mathbf{a}}\in\mathcal{D}_{[N]}([-1,1])$ such that $\mathcal{B}(f_{\mathbf{a}})=\mathbf{a}$ and for all flat $g\in\mathcal{D}_{[N]}([-1,1])$ we get:
$$\mathcal{B}(g+f_{\mathbf{a}})=(g^{(k)}(0)+f_{\mathbf{a}}^{(k)}(0))_{k\in\NN}=\mathbf{a}.$$
We consider $\varphi_c$ with $c\ge 3a$ and set
$$\psi_c(x):=\varphi_c\left(x-\frac{1}{2}\right),\;\;\; x\in[-\frac{1}{3},\frac{1}{3}],\hspace{20pt}\psi_c(x)=0,\;\;\;x\in[-1,1]\backslash[-\frac{1}{3},\frac{1}{3}],$$
hence $\psi_c\in\mathcal{D}_{[N]}([-1,1])$ with $\supp(\psi_c)=[-\frac{a}{c},\frac{a}{c}]$ and each $\psi_c$ is flat at $x=0$. We define the linear span:
$$\mathcal{V}:=\left\{\sum_{i=1}^l\alpha_i\psi_{c_i}: l\in\NN_{>0}, \alpha_i\in\CC, 3a<c_1<\dots<c_l\right\},$$
which has all the required properties since the functions $\psi_{c_i}$ are linearly independent.
\qed\enddemo

Using $\varphi_c$ from Lemma \ref{Petzsche2} we can also prove a multiplicative variant of the previous statement. Here, for given $\mathbf{a}=(a_k)_{k\in\NN}$ and any $c>0$ we note that $\mathcal{B}(\varphi_c\cdot f_{\mathbf{a}})=\mathbf{a}$ because
$$\forall\;k\in\NN:\;\;\;(\varphi_c\cdot f_{\mathbf{a}})^{(k)}(0)=\sum_{0\le j\le k}\binom{k}{j}\underbrace{\varphi_c^{(j)}(0)}_{=\delta_{j,0}}\cdot f_{\mathbf{a}}^{(k-j)}(0)=f_{\mathbf{a}}^{(k)}(0)=a_k.$$
Fix now $c\ge a$. Then $\varphi_c\cdot f_{\mathbf{a}}\in\mathcal{D}_{[N]}([-1,1])$. Indeed, here $b\le 1$ and since $\varphi_c,f_{\mathbf{a}}\in\mathcal{D}_{[N]}([-1,1])$ the log-convexity for $N$ implies stability under pointwise multiplication for this class, e.g. see the proof of \cite[Prop. 2.0.8]{diploma}.

Then define
$$\mathcal{V}_1:=\left\{\varphi_{c}+\sum_{i=1}^l\lambda_i\psi_{c_i}: \lambda_i\in\CC, l\in\NN_{>0}, 3a<c_1<\dots<c_l\right\}$$
and note that $\mathcal{V}_1$ is only an affine space. When considering $\mathbf{a}=0$, then we can construct an infinite-dimensional vector space as follows: We can find $f\in\mathcal{D}_{[N]}([-1,1])$ such that $\mathcal{B}(f)=\mathbf{a}$. Since $N$ is nonquasianalytic, we can assume that $f\neq 0$ (e.g. take any function with compact support contained in $[-1,1]$ and such that $0\notin\supp(f)$). Then the space
$$\mathcal{V}_2:=\left\{\sum_{i=1}^l\lambda_i\varphi_{c_i}: \lambda_i\in\CC, l\in\NN_{>0}, a<c_1<\dots<c_l\right\}$$
has all the required properties.

\begin{remark}\label{noninjremark}
For completeness we mention that (the proof of) Proposition \ref{noninj} can be transferred to the $r$-ramified setting ($r\in\NN_{>0}$) as well by taking functions $\varphi_c\in\mathcal{D}_{r,[N]}([-1,1])$ (when this space is nontrivial) with $\varphi_c^{(j)}(0)=\delta_{j,0}$. The existence of such functions is ensured by Lemma \ref{Petzsche2} applied to the so-called $r$-interpolating sequence $P^{N,r}$, see \cite[Sect. 2.5, Lemma 3.5]{mixedramisurj} and \cite[Lemma 2.3]{Schmetsvaldivia00}.\vspace{6pt}

However, the ''multiplicative variant'' seems to be unclear since in general $\mathcal{D}_{r,[N]}([-1,1])$ is not closed under the pointwise product of functions.
\end{remark}

\subsection{Failure of surjectivity in the nonquasianalytic setting}\label{failuresurj}
Our aim is to show that in the situation of Theorem \ref{SVtheorem} the inclusion $\mathcal{B}(\mathcal{D}_{\{N\}}([-1,1]))\supseteq\Lambda_{\{M\}}$ is strict for all $N$ which are not strongly nonquasianalytic. More precisely, the complement $\mathcal{B}(\mathcal{D}_{\{N\}}([-1,1]))\backslash\Lambda_{\{M\}}$ is large in the following sense, see \cite[Def. 1.4]{BPS} and the references therein:\vspace{6pt}

A set $\mathcal{L}$ in a vector space $V$ is called {\itshape lineable} in $V$ if $\mathcal{L}\cup\{0\}$ contains an infinite-dimensional vector space.\vspace{6pt}

We start with the following preparation: Let $N\in\hyperlink{LCset}{\mathcal{LC}}$. Then the function $\theta_N$ defined by
$$\theta_{N}(x):=\sum_{k=0}^{+\infty}\frac{N_k}{(2\nu_k)^k}\exp(2i\nu_k
x), \quad x \in \RR,$$
satisfies the following conditions: $\theta_N\in\mathcal{E}_{\{N\}}(\RR,\CC)$ and
\begin{equation}\label{thetafunction}
\theta_N^{(j)}(0)=i^j s_{j}  \text{ with } s^1_j:=s_{j}\ge N_{j}\, ,  \quad \forall\,j\in\NN.
\end{equation}
We refer to \cite[Theorem 1]{thilliez}; for a detailed proof see also \cite[Prop. 3.1.2]{diploma} and \cite[Lemma 2.9]{compositionpaper}. It is not difficult to see that $\theta_N$ does not belong to the Beurling type class $\mathcal{E}_{(N)}(\RR,\CC)$.

Using this function we can show the following:

\begin{proposition}\label{nqfailure}
Let $N\in\hyperlink{LCset}{\mathcal{LC}}$ such that \hyperlink{gamma1}{$(\gamma_1)$} fails. Then for any sequence $M\in\mathcal{N}_{\preceq,R}$ with $\mathcal{B}(\mathcal{D}_{\{N\}}([-1,1]))\supseteq\Lambda_{\{M\}}$ we have strict inclusion, more precisely the set
$$\mathcal{B}(\mathcal{D}_{\{N\}}([-1,1]))\backslash\bigcup_{M\in\mathcal{N}_{\preceq,R}}\Lambda_{\{M\}}$$
is lineable in $\mathcal{B}(\mathcal{D}_{\{N\}}([-1,1]))$ (the constructed vector space has dimension $\aleph_0$).


\end{proposition}

\demo{Proof}
By Theorem \ref{SVtheorem}, for any $M\in\mathcal{N}_{\preceq,R}$ with $\mathcal{B}(\mathcal{D}_{\{N\}}([-1,1]))\supseteq\Lambda_{\{M\}}$ condition $(M,N)_{SV}$ is satisfied, hence $N$ is nonquasianalytic and by Lemma \ref{optimalSV} in order to conclude it is sufficient to show lineability of $\mathcal{B}(\mathcal{D}_{\{N\}}([-1,1]))\backslash\Lambda_{\{L\}}$ for the maximal sequence $L\equiv L^1$ (see \eqref{optimalSVsequence}).\vspace{6pt}

Take $\theta_N\in\mathcal{E}_{\{N\}}(\RR,\CC)$ and $\psi\in\mathcal{D}_{\{N\}}([-1,1])$ such that $\psi^{(j)}(0)=\delta_{j,0}$ and $\supp(\psi)\subseteq[-1,1]$ (e.g. take $\psi\equiv\varphi_c$, $c\ge a:=\sum_{j=1}^{+\infty}\frac{1}{\nu_j}$) and set
$$\Theta_N:=\theta_N\cdot\psi.$$
Then $\Theta_N\in\mathcal{D}_{\{N\}}([-1,1])$ because log-convexity of $N$ implies that the class is closed under pointwise multiplication (see e.g. the proof of \cite[Prop. 2.0.8]{diploma}). Moreover,
$$\forall\;k\in\NN:\;\;\;\Theta_N^{(k)}(0)=\sum_{0\le j\le k}\binom{k}{j}\theta_N^{(j)}(0)\psi^{(k-j)}(0)=\theta_N^{(k)}(0)=i^ks_k.$$
\vspace{6pt}

{\itshape Claim I:} $\mathcal{B}(\Theta_N)\in\mathcal{B}(\mathcal{D}_{\{N\}}([-1,1]))\backslash\Lambda_{\{L\}}$.

Clearly $\mathcal{B}(\Theta_N)\in\mathcal{B}(\mathcal{D}_{\{N\}}([-1,1]))$. Since $N$ fails \hyperlink{gamma1}{$(\gamma_1)$}, by Lemma \ref{gamma1} we know that $L\hyperlink{approx}{\approx}N$ is violated. Now assume that there exists $\mathbf{a}=(a_j)_j\in\Lambda_{\{L\}}$ with $\mathcal{B}(\Theta_N)=\mathbf{a}$. Then
$$\exists\;C,h>0\;\forall\;j\in\NN:\;\;\;N_j\le|\Theta_N^{(j)}(0)|=|a_j|\le C h^jL_j,$$
which implies $N\hyperlink{preceq}{\preceq}L$. By Lemma \ref{optimalSV} this implies $L\hyperlink{approx}{\approx}N$, a contradiction.\vspace{6pt}

{\itshape Claim II:} $\mathcal{B}(\mathcal{D}_{\{N\}}([-1,1]))\backslash\Lambda_{\{L\}}$ is lineable in $\mathcal{B}(\mathcal{D}_{\{N\}}([-1,1]))$.

First, for each $h>0$ we set $N^h:=(h^jN_j)_{j\in\NN}$. Obviously $N^h\hyperlink{approx}{\approx}N$, $\nu^h_j:=\frac{N^h_j}{N^h_{j-1}}=h\nu_j$, $j\in\NN_{>0}$, and $\nu^h_0:=1$. Moreover each $N^h$ is clearly log-convex and normalized for all $h\ge 1$. So it makes sense to consider $\theta_{N^h}$ and $\Theta_{N^h}:=\theta_{N^h}\cdot\psi\in\mathcal{D}_{\{N^h\}}([-1,1])=\mathcal{D}_{\{N\}}([-1,1])$. By the previous comments we get
$$\forall\;h\ge 1\;\exists\;A_h,B_h\ge 1\;\forall\;j\in\NN:\;\;\;h^jN_j=N^h_j\le|\Theta_{N^h}^{(j)}(0)|\le A_hB_h^jN^h_j=A_h(B_hh)^jN_j,$$
and
$$\forall\;j\in\NN\;\forall\;h\ge 1:\;\;\;\Theta_{N^h}^{(j)}(0)=i^js_j^h,\;\;\;s_j^h\ge N^h_j=h^jN_j.$$
We introduce iteratively a sequence of functions $(\Phi_j)_{j\ge 1}$ and strictly increasing sequences of numbers $(A_j)_{j\ge 1}$, $(B_j)_{j\ge 1}$, $(C_j)_{j\ge 1}$ such that $B_{j+1}>C_j>B_j(>\max_{1\le i\le j-1}B_i$) as follows: First set
$$\Phi_1:=\Theta_N(=\Theta_{N^1}),$$
which implies
$$\exists\;A_1,B_1\ge 1\;\forall\;l\in\NN:\;\;\;N_l\le|\Phi_1^{(l)}(0)|=|\Theta_N^{(l)}(0)|\le A_1B_1^lN_l.$$
Then we put iteratively
$$\Phi_{j+1}:=A_j\cdot\Theta_{N^{C_j}},\;j\in\NN_{>0},$$
satisfying
$$\exists\;A_{j+1},B_{j+1}\ge 1\;\forall\;l\in\NN:\;\;\;A_jN_l^{C_j}=A_jC_j^lN_l\le|\Phi^{(l)}_{j+1}(0)|\le A_{j+1}B_{j+1}^lN_l.$$
The choices of the sequences yield
\begin{equation*}\label{Phiequ}
\forall\;l\in\NN\;\forall\;j\in\NN_{>0}:\;\;\;\Phi^{(l)}_{j+1}(0)=A_j\Theta_{N^{C_j}}^{(l)}(0)=A_ji^ls_l^{C_j},\;\;\;s_l^{C_j}\ge N_l^{C_j}=C_j^lN_l,
\end{equation*}
and $\Phi_1^{(l)}(0)=i^ls_l=i^ls^1_l$, $s^1_l\ge N_l$. 
Finally, we define
$$\mathcal{V}:=\left\{\sum_{i=1}^k\alpha_i\Phi_i, k\in\NN_{>0}, \alpha_i\in\CC\right\},$$
so $\mathcal{V}\subseteq\mathcal{D}_{\{N\}}([-1,1])$ and the functions $\Phi_i$ are linearly independent. It is then straightforward that $\mathcal{B}(\mathcal{V})\cap\Lambda_{\{L\}}=\{0\}$.

\qed\enddemo

\bibliographystyle{plain}
\bibliography{Bibliography}
\end{document}